\documentclass[twoside,a4paper,12pt,centertags]{amsart}

\usepackage[T1]{fontenc}
\usepackage[utf8]{inputenc}
\usepackage{subfiles}
\usepackage{multicol}

\makeatletter
\usepackage[section]{placeins}

\AtBeginDocument{%
	\expandafter\renewcommand\expandafter\subsection\expandafter{%
		\expandafter\@fb@secFB\subsection
	}%
}
\makeatother
\makeatletter
\@namedef{subjclassname@2020}{%
	\textup{2020} Mathematics Subject Classification}
\makeatother

\usepackage{amsfonts,amssymb,amsmath,amsthm}
\usepackage{mathrsfs}
\usepackage{nicefrac}
\usepackage{latexsym}
\usepackage{mathtools}
\usepackage[noadjust]{cite}
\usepackage{paralist}
\usepackage{longtable}
\usepackage{float}
\usepackage{pdfsync}
\usepackage{tikz,tikz-cd}
\usetikzlibrary{patterns}
\usepackage{hhline}
\usepackage{import}
\usepackage{esint}

\usepackage{enumitem}

\usepackage[colorlinks=true,citecolor=magenta,linkcolor=cyan,backref]{hyperref}
\AtBeginDocument{}

\numberwithin{equation}{section}

\theoremstyle{plain}
\newtheorem{thm}{Theorem}[section]
\newtheorem{lem}[thm]{Lemma}
\newtheorem{cor}[thm]{Corollary}
\newtheorem{prop}[thm]{Proposition}

\theoremstyle{definition}
\newtheorem{defn}[thm]{Definition}
\newtheorem{rem}[thm]{Remark}

\usepackage{color}

\DeclareMathOperator*{\PV}{PV}
\renewcommand{\div}{\operatorname{div}}

\newcommand{\IC}{\mathbb{C}}

\newcommand{\R}{\mathbb{R}}

\newcommand{\cA}{\mathcal{A}}
\newcommand{\cT}{\mathcal{T}}
\newcommand{\cD}{\mathcal{D}} 

\newcommand{\cK}{\mathcal{K}} 
\newcommand{\cL}{\mathcal{L}}
\newcommand{\cM}{\mathcal{M}}

\newcommand{\cP}{\mathcal{P}}

\newcommand{\cS}{\mathcal{S}} 

\newcommand{\cB}{\mathcal{B}}

\newcommand{\loc}{\operatorname{loc}}
\renewcommand{\L}{\operatorname{L}} 
\newcommand{\C}{\operatorname{C}} 
\renewcommand{\H}{\operatorname{H}} 
\newcommand{\B}{\operatorname{B}}
\newcommand{\F}{\operatorname{F}}
\newcommand{\E}{\operatorname{E}}%
\newcommand{\G}{\operatorname{G}}

 \newcommand{\X}{\operatorname{X}} 
 \newcommand{\Y}{\operatorname{Y}}
  \newcommand{\I}{\operatorname{I}}

   \newcommand{\Z}{\operatorname{Z}}
    \newcommand{\p}{\operatorname{P}}

\DeclareRobustCommand{\Hdot}{\dot{\H}\protect{\vphantom{H}}
} 
\DeclareRobustCommand{\Edot}{\dot{\E}\protect{\vphantom{E}}} 
\DeclareRobustCommand{\Fdot}{\dot{\F}\protect{\vphantom{F}}} 
\DeclareRobustCommand{\Gdot}{\dot{\G}\protect{\vphantom{G}}} 
\DeclareRobustCommand{\Bdot}{\dot{\B}\protect{\vphantom{B}}} 
 
 \DeclareRobustCommand{\Pdot}{\dot{\p}\protect{\vphantom{P}}}

\DeclareRobustCommand{\Ydot}{\dot{\Y}\protect{\vphantom{Y}}} 
\DeclareRobustCommand{\Zdot}{\dot{\Z}\protect{\vphantom{Z}}}

\newcommand{\Bd}{\Bdot}

\newcommand{\Fb}{\F^{\beta}}
\newcommand{\Fbd}{\Fdot^{\beta}}
\newcommand{\Gb}{\G^{\beta}}
\newcommand{\Gbd}{\Gdot^{\beta}}
\newcommand{\Hb}{\H^{\beta}}
\newcommand{\Hbd}{\Hdot^{\beta}}
\newcommand{\Vb}{\Y^{\beta}}
\newcommand{\Vbd}{\Ydot^{\beta}}
\newcommand{\Vbk}{\Y^{\beta}_{\ka}}

\newcommand{\Wb}{\Z^{\beta}}
\newcommand{\Wbd}{\Zdot^{\beta}}
\newcommand{\Wbk}{\Z^{\beta}_{\ka}} 

\newcommand{\Fbk}{\F^{\beta}_{\ka}}
\newcommand{\Bka}{\B_{\ka}}

\renewcommand\Re{\operatorname{Re}}

\def\angle#1#2{\langle #1,#2 \rangle} 
\newcommand{\llangle}{\langle\!\langle}
\newcommand{\rrangle}{\rangle\!\rangle}

\def\tpsi{\tilde\psi}

\newcommand\ka{\kappa}
\newcommand\gb{\gamma({\beta})}
\newcommand\essinf{\operatorname{ess\, inf}}

\renewcommand{\iint}{\int_{}\kern-.34em \int} 
\renewcommand{\iiint}{\iint_{}\kern-.34em \int} 

\newcommand{\dd}{\, \mathrm{d} }
\newcommand{\dv}{\, \mathrm{d} v}
\newcommand{\dx}{\, \mathrm{d} x}
\newcommand{\ds}{\, \mathrm{d} s}

\newcommand{\dt}{\, \mathrm{d} t}

\def\oa{\mathbf{a}}
\def\ob{\mathbf{b}}
\def\oc{\mathbf{c}}

\title[Fundamental solutions to Kolmogorov-Fokker-Planck equations]{Fundamental solutions to Kolmogorov-Fokker-Planck equations with  unbounded lower order terms}
\author{Pascal Auscher}
\address{Universit\'e Paris-Saclay, CNRS, Laboratoire de Math\'{e}matiques d'Orsay, 91405 Orsay, France}
\email{pascal.auscher@universite-paris-saclay.fr}
\date{November 7, 2026}

\thanks{}



\begin{document}

\begin{abstract} We continue the study  of fundamental solutions for equations involving Kolmogorov-Fokker-Planck operators with rough local or non-local diffusions. Here we consider possibly unbounded lower order terms, such as drifts and potentials. We prove existence and uniqueness of  weak solutions and  of an appropriate notion of fundamental solution for a wide class of kinetic equations.
This is based on kinetic transfer of regularity and anisotropic embeddings.

\end{abstract}

\maketitle

\setcounter{tocdepth}{1}
\tableofcontents

\section{Introduction}

Consider a local elliptic operator
\[ 
\cA= -\div_{v}A \nabla_{v} -\div_{v}\oa  +\ob\cdot  \nabla_{v}+\oc
\]
where  $A$, $\oa$, $\ob$, $\oc$ depend  on time $t\in \R$, position $x\in \R^d$, velocity $v\in \R^d$,  $A$ being a measurable, bounded and elliptic matrix. We ask under what conditions on the lower order terms $\oa$, $\ob$, $\oc$, the rough Fokker-Planck operator
\[\cL=(\partial_{t}+v\cdot\nabla_{x})+\cA\] 
has a fundamental solution.  Such operators are models for the linearization of the Landau equation. 

One can ask the same question starting from a non-local diffusive operator such as
\begin{equation}
  \label{eq:A-kernel}
  \cA_{\beta} f (t,x,v) = \PV \int_{\R^d} \frac{f(t,x,v)  - f(t,x,v')}{|v-v'|^{d+2\beta}}  A_{\beta} (t,x,v,v') \dv',  \quad 0<\beta<1,
\end{equation}
where $A_{\beta}$ is bounded, measurable and uniformly coercive. What kind of lower order perturbation  $\cP$ (potential, lower order fractional operators) can we add and still define a fundamental solution for the Kolmogorov operator  $(\partial_{t}+v\cdot\nabla_{x})+\cA_{\beta}+\cP$?
Such operators are models for  the linearization of the Boltzmann equation. 

\

The usual approach is to look for existence of a kernel function such that the integral formula 
\begin{equation}
\label{eq:representation}
 f(t,x,v)=\int_{(-\infty,t)\times \, \R^{2d}} \Gamma(t,x,v,s,y,w)S(s,y,w)\dd(s,y,w)
\end{equation}
 for good sources $S$ produces (weak or strong) solutions to $((\partial_{t}+v\cdot\nabla_{x})+\cA)f=S$ and then to prove estimates and establish  properties for this kernel, customarily called fundamental solution in the literature coming from the terminology of PDEs with constant coefficients (and we shall use the terminology fundamental solution kernel  to make the distinction with the notion we use here). Another question is uniqueness of this object and also under which condition on a solution $f$ can one write the integral formula.  For local diffusions, many works exist on these topics but the question has not been fully  resoleved when  perturbations with unbounded coefficients are added to the higher order term. The situation for non-local diffusions has just been answered recently without perturbations. See below. 
 
 \

Similarly to parabolic equations, the topic of fundamental solutions is often inextricably associated to regularity  theory  for solutions (pointwise estimates, Harnack inequalities, Moser estimates, De Giorgi methods, Nash estimates...). In the important works of Aronson \cite{Ar67}, \cite{Ar68}, the construction of  generalized fundamental solutions for parabolic equations relies on Nash's result. Thinking that    \eqref{eq:representation} can be formally written as the operational  equality
\begin{equation}
\label{eq:representation2}
f(t,\cdot)=\int_{-\infty}^t\Gamma(t,s)S(s,\cdot)\, \dd s,
\end{equation}
the operator $\Gamma(t,s)$ should be the operator mapping data at time $s$ to  the solution at time $t$ of $((\partial_{t}+v\cdot\nabla_{x})+\cA)f=0$ on $(s,\infty)\times \, \R^{2d}$.
One point we want to emphasize is that fine regularity theory is not necessary to construct $\Gamma(t,s)$ and to obtain the  representation \eqref{eq:representation2}  substituting \eqref{eq:representation}, see \cite{AB} for implementation in the parabolic setting inspired by and extending the early works of Lions \cite{Lions}.  Regularity theory under additional assumptions intervenes afterwards to obtain further properties of  that fundamental solution, so as to reach the kernel representation \eqref{eq:representation} with estimates.  Uniqueness guarantees that any other method, for instance through regularisation of coefficients, necessarily leads to the kernel of $\Gamma (t,s)$. 

\

Historically, for local  diffusions having smooth coefficients (for example, H\"older-continuous coefficients), one can construct the fundamental solution kernel by Levi's parametrix method and prove estimates. See the references in the introduction of \cite{MR2000352}.
In that article,  pointwise estimates are established  assuming H\"older-continuous, real, symmetric leading coefficients (with no lower order terms),  following Aronson's method.  An a priori bound on the fundamental solution kernel not using any modulus of continuity is obtained, but existence of the fundamental solution kernel not assuming regularity was not proved. This was extended to  lower order perturbations with bounded coefficients in \cite{MR3717350}, 
and more recently in \cite{francesca2021fundamental} 
with some unbounded coefficients (for example $\oc \in \L^{q}_{t,x,v}$ when $q>2d+1$). Existence of a fundamental solution kernel was proposed  in this article under divergence free structural restrictions and real symmetry on the leading coefficients.     The situation was resolved only recently in \cite{AIN2} for local diffusions with  bounded lower order coefficients, following the principles described above. First, only under coercivity assumptions, not even assuming real-valued coefficients, there  exists a unique fundamental solution in the sense of a family of operators $\Gamma(t,s)$.  Secondly, the existence of a kernel function to such fundamental solution operators together with upper estimates follows from (is in fact shown to be equivalent to)  local scale invariant $\L^2-\L^\infty$ Moser estimates for the local weak solutions of $\cL f=0$ and those of its backward adjoint. When the leading coefficients are real and symmetric with no lower order terms, these Moser estimates were obtained in \cite{MR2068847} with a stronger notion of solutions that fits smooth coefficients. 
See also \cite{MR3923847}  with a notion of weak solution that fits what we use here. Symmetry is not needed in these arguments. 
The conclusion in the case of Fokker-Planck equations with only leading real measurable coefficients without  lower order bounded coefficients admitting Moser estimates  is that there exists a fundamental solution kernel. And it has the expected sharp pointwise upper estimates in terms of the kinetic geometry. 
Adding lower order terms with bounded coefficients gives the same conclusion provided the Moser estimates are established.  Related results with different methods are in \cite{Gaddy}. 

\

For non-local diffusions, the situation is much less understood. This was only recently  done  
  for the Kolmogorov operators  $(\partial_{t}+v\cdot\nabla_{x})+ (-\Delta_{v})^{\beta}$,  $0<\beta<1$. While existence of the fundamental solution kernel (by Fourier methods or probabilistic methods) is not an issue, the sharp upper  pointwise decay with a polynomial profile taking into account transport and sub-diffusion in velocity    was obtained  in \cite{Hou2026}. This corresponds to  $A_{\beta}$ constant in \eqref{eq:A-kernel} and even $A_{\beta}=\kappa(v-v')$ is considered for a class of positive bounded functions $\kappa$.  For more general non-local diffusions,  there is an actual quest for understanding the regularity theory. For example, Harnack inequalities can fail \cite{MR4049224}
   and should be replaced by inequalities involving non-local behavior (see \cite{anceschi2026giorginashmosertheorykineticequations}). 
 Concerning fundamental solutions,  the same abstract strategy as for the local diffusions applies and yields the existence of  fundamental solution operators $\Gamma(t,s)$ for a wider class of non-local diffusions than just those defined by \eqref{eq:A-kernel}, see  \cite{AIN2}. And in  collaboration with L. Niebel \cite{AN2}, under the assumption $0<\lambda\le A_{\beta}(t,x,v,v')=A_{\beta}(t,x,v',v)\le \Lambda$ for almost every $(t,x,v,v')$, we obtained the existence of a kernel that is nonnegative and enjoys  the same sharp   pointwise upper estimate as in \cite{Hou2026}. 
  
  \
  
  In this work, we attack the problem of existence of  fundamental solution operators $\Gamma(t,s)$ for kinetic operators with both local and non-local diffusions  under perturbation by lower order terms having unbounded coefficients (potentials, drifts and related terms). For instance, the size of the coefficients can be  measured in mixed Lebesgue-Lorentz spaces in the variables $(t,x,v)$. This is obtained via an abstract strategy based on  transfer of regularity and embeddings of \cite{AIN2} to handle the lower order perturbations in a short and systematic way. Despite similarities, this is not an extension of \cite{AB} as the initial input of a variational space does not apply here. \footnote{The article \cite{KN} posted one day after ours proved such a variational inequality (see Prop. 3.2 there).  Still, to implement the strategy of \cite{AB} in our applications would require  to establish  embeddings that are more complicated  due to non commutativity of the kinetic group law.}

  Our strategy  applies to many examples with critical scaling, such as Coulomb potentials in velocity and also potential with spatial decay  reflecting the transfer of regularity. For instance, for local diffusions, one can take potentials $\oc$ on the order of the Coulomb potential $c|v|^{-2}$ with $ \Re c$ bounded below from the Hardy inequality 	and $|c|$ bounded above, or the same with an additional $c'|x|^{-2/3}$ perturbation,   $|c'|$ small,  or $\oc\in  \L^{2d+1}_{t,x,v}$ of arbitrary size, etc; for non-local diffusion:  $\oc$ on the order of  $c|v|^{-2\beta}$ with $ \Re c$ bounded below from  the Hardy-Rellich inequality and $|c|$ bounded above, the same with an additional   $ c'|x|^{-2\beta/(2\beta+1)}$ with $|c'|$ small or $\oc\in  \L^{(1+1/\beta)d+1}_{t,x,v}$ of arbitrary size, etc.
  
 \
 
 The results also apply to parabolic equations with no position variable. In the kinetic setting, the transfer of regularity that takes into account the transport flow is an essential tool.  
  
  The question of existence of a representing kernel with estimates for kinetic equations will be treated in a subsequent work.

 \subsection*{Acknowledgments}  The author acknowledges discussions with C. Imbert, and with L. Niebel, who also helped for the references. ChatGPT-6Astra was used only to find the references for the sharp constant in the Hardy-Rellich inequality, and to eliminate mathematical and text  typos, and grammatical errors.  
 
A CC BY-NC 4.0 \url{https://creativecommons.org/licenses/by/4.0/} public copyright license has been applied by the authors to the present document and will be applied to all subsequent versions up to the Author Accepted Manuscript arising from this submission.

 \section{Functional Spaces}

 We work on the full space $\R\times \,\R^{2d}$, half-spaces $(0,\infty)\times \,\R^{2d}$ or strips $(0,T)\times \,\R^{2d}$. If $\I$ is an open interval, we call $\C^{}_{0}\!\X= \C^{}_{0}(\overline \I \, ; \X)$,  the space of continuous functions on $\overline \I$ into $\X$ with 0 limit at the infinite ends of $\I$. We do not indicate the time interval $\I$ unless it is necessary for the understanding.  The Lebesgue spaces are understood with respect to  Lebesgue measure. We indicate the variables in the notation. 
 The suffix $'$ means the dual space for  the duality extending the $\L^2_{t,x,v}$ inner product. 

 The $\L^2$ based Sobolev spaces  of order $\alpha\in \R$ on $\R^d$ are denoted by $\H^{\alpha}$. Throughout this article  $\beta>0$ and no other assumptions are made unless specified.    Define
\begin{align*}
\Fb&= \L^2_{t,x}\H^{\beta}_{\vphantom{t^a} v},
\\
\Gb&= \L^2_{t,v}\H_{\vphantom{t} x}^{\gb}, \quad  \gb=\tfrac{\beta}{2\beta+1}\\
\Vb&= \Fb\cap \Gb\cap \C^{}_{0}\L^2_{x,v}
\\
 \Wb&=(\Fb)'+  (\Gb)'+ \L^1_{t} \L^2_{x,v} 
\end{align*}
equipped with the natural norms.    We use the space $\cS_{K}$ defined as  the subset of those $S\in \cD_{t}\, \cS_{x,v}$ whose partial Fourier transform  has compact support in $\I^{}_{t}\times \, (\R^{d}_{\varphi}\setminus \{0\})\times\,  (\R^{d}_{\xi}\setminus \{0\})$. This space is dense in all the spaces listed above.  Note that $\Vb$ is a closed subspace of  $(\Wb)'=\Fb\cap \Gb\cap \L^\infty_{t} \L^2_{x,v} $.

We shall use the homogeneous  Sobolev spaces $\Hdot^\alpha$ on $\R^d$ as defined in \cite{AN}, to which we refer for details,  via the Littlewood-Paley decomposition  as a subset of the space of tempered distributions $\cS'$.  
On the dense subset of Schwartz functions with compact support in $\R^d\setminus\{0\}$, the norm is comparable to  $\||\xi|^{\alpha}\hat f\|_{\L^2}$ via  Fourier transform. When $\alpha<d/2$, all elements of   $\Hdot^\alpha$ have Fourier transform equal to the product of an $\L^2$ function with $|\xi|^{-\alpha}$.  Define the homogeneous versions of the previous spaces
\begin{align*}
\Fbd&= \L^2_{t,x}\Hbd_{\vphantom{t^a} v},
\\
\Gbd&= \L^2_{t,v}\Hdot_{\vphantom{t} x}^{\gb}\\
\Vbd&= \Fbd\cap \Gbd\cap \C^{}_{0}\L^2_{x,v} 
\\
 \Wbd&=(\Fbd)'+  (\Gbd)'+ \L^1_{t} \L^2_{x,v} 
\end{align*}
equipped with the natural norms (semi-norm for $\Fbd$ if $\beta\ge d/2$).   The space $\cS_{K}$ is dense in all of the listed spaces.   
Note that $\Vbd$ is a closed subspace of  $(\Wbd)'=\Fbd\cap \Gbd\cap \L^\infty_{t} \L^2_{x,v} $. There is a subtlety in the definition of $\Fbd$ which implies that, if $\beta\ge d/2$, polynomials $P=P(v)$ in the $v$-variable cannot belong to $\Fbd$ because this space is defined as a subspace of $\L^2_{t,x}\cS'_{v}$. For $\Gbd$, this is always a normed space as $\gb<1/2$. 

The associated (semi-)norms are homogeneous with respect to the kinetic scaling $h_{\delta}:(t,x,v) \mapsto (\delta^{2\beta} t, \delta^{2\beta+1} x, \delta v)$ for any  
$\delta>0$, with  $$\|f\circ h_{\delta}\|= \delta^{-(\beta+1)d}\|f\|$$ 
for the $\Fbd, \Gbd, \Vbd$ norms and
$$\|f\circ h_{\delta}\|= \delta^{-2\beta-(\beta+1)d}\|f\|$$ 
for the $(\Fbd)',  (\Gbd)', \Wbd$ norms.

The inhomogeneous versions  $\Fb, \Gb, \Vb$ are equal to the homogeneous versions  $\Fbd, \Gbd, \Vbd$ intersected with $\L^2_{t,x,v}$,  adding to the homogeneous norm a multiple factor of the $\L^2_{t,x,v}$ norm. It will be useful to quantify this factor, appearing as a spectral parameter $\ka$ in the estimates.  We shall introduce the notation after stating our main hypotheses. 

For $(\Fb)', (\Gb)', \Wb$, we add the $\L^2_{t,x,v}$ space to  $(\Fbd)', (\Gbd)', \Wbd$, hence they  are larger spaces.

 \section{Review on kinetic embedding and transfer of regularity}
 
Embeddings and transfer of regularity are the main \textit{a priori} inequalities used to obtain existence and uniqueness results for kinetic equations. In the case of perturbations with lower order terms having unbounded coefficients on infinite strips, we must also track  the dependence in a spectral parameter that measures inhomogeneity, which we denote by $\ka$. This was not crucial when the perturbations had bounded coefficients.

We do not indicate the open interval to which the time variable $t$ belong in the subsequent notation. It is implicit.

\begin{thm}[Spectrally quantified kinetic embedding and transfer of regularity]\label{thm:optimalregularity} 
  Let $\beta>0$. Let $\I$ be an open interval in $\R$.  Let $\ka\in [0,\infty)$. Assume $\ka$ is arbitrary when $\I$ is unbounded and $\ka^2|\I | \sim 1$ when $\I$ is bounded. 
  Assume 
  \begin{align*}
 f &\in \L^2_{t,x}\Hdot^{\beta}_{\vphantom{t,x} v}, \quad  \ka f\in \L^2_{t,x,v}
 \\
 (\partial_{t}+v\cdot\nabla_{x})f &=S_{1}+S_{2}+S_{3}+ \ka S_{4}
 \\
 S_1 \in \L^2_{t,x}\Hdot^{-\beta}_{\vphantom{t,x} v}, \quad &S_2 \in \L^2_{t,v}\Hdot_{\vphantom{t,x} x}^{-\gb}, \quad S_3 \in \L^1_{\vphantom{t,x} t}\L^2_{\vphantom{t,x} x,v}, \quad S_{4}\in \L^2_{\vphantom{t,x} t,x,v}.
\end{align*}
Then, we have
\begin{enumerate}
\item $f\in  \L^2_{t,v}\Hdot_{\vphantom{t,v} x}^{\gb} \cap \C_{0}(\overline \I\, ; \, \L^2_{x,v})$ with 
\begin{align} \label{eq:embedding}
\| f \|_{\L^2_{t,v}\Hdot_{\vphantom{t,v} x}^{\gb}}+ &\sup_{t\in \overline{\I}} \|f(t)\|_{\L^2_{x,v}} 
\lesssim_{d,\beta} \|f\|_{\L^2_{t,x} \Hdot^{\beta}_{\vphantom{t,x} v}}
+  \ka \|f\|_{\L^2_{\vphantom{t,x} t,x,v}}
\\ \nonumber
&+ \|S_1 \|_{\L^2_{t,x}\Hdot^{-\beta}_{\vphantom{t,x} v}}
+ \|S_2 \|_{\L^2_{t,v}\Hdot_{\vphantom{t,x} x}^{-\gb}}
+ \|S_3 \|_{\L^1_{\vphantom{t,x} t} \L^2_{x,v}} 
+  \|S_{4}\|_{\L^2_{\vphantom{t,x} t,x,v}}.
\end{align}
\item The map  $t\mapsto\|f(t)\|_{\L^2_{x,v}}^2$ is absolutely continuous on $\I$ with for a.e. $t\in \I$, 
\begin{equation}\label{e:abs-cont}
\frac{\mathrm{d} }{\mathrm{d}t }\|f(t)\|_{\L^2_{x,v}}^2 =   2\Re \bigg(\int_{\R^d} \angle {S_{1}}{f} \dx + \int_{\R^d} \angle {S_{2}}{f} \dv  + \iint_{\R^{2d}}  {S_{3}}\overline f +  \ka {S_{4}}\overline f  \, \dx\dv \bigg).
\end{equation}
\end{enumerate}
In \eqref{e:abs-cont}, the first bracket denotes the sesquilinear duality between $\Hdot^{-\beta}_{\vphantom{t}v}$ and $\Hdot^{\beta}_{\vphantom{t}v}$, giving an integrable function of $(t,x)$,  and the second  one refers to the duality between  $\Hdot_{\vphantom{t}x}^{-\gb}$ and $\Hdot_{\vphantom{t}x}^{\gb}$, giving an integrable function of $(t,v)$. Integrals are in the sense of Lebesgue.
\end{thm}

It is important to note that the implicit constant in \eqref{eq:embedding} depends on $d$ and $\beta$ but not on $f$, the interval $\I$ and the spectral parameter $\ka$, especially when $\I$ is unbounded. The case $\ka=0$ corresponds to a homogeneous situation with $S_{4}=0$.

\begin{proof} 
Via Fourier transform in the $x,v$ variables, we define the Kolmogorov integral operators
 $$\cK^\pm_{\beta,\ka}= \big(\pm (\partial_{t}+v\cdot\nabla_{x}) + (-\Delta_{v})^{\beta}+\ka^{2} \, \big)^{-1}, \quad  \ka\ge 0.$$
 
 \

{\em Step 1:  $\ka=0$.} 

  For $\I=\R$, the statement is that of  \cite[Theorem~1.16]{AIN} when $\beta<d/2$. The supplementary assumption when $\beta\ge d/2$ was removed in  \cite[Theorem~7.1]{AN}, see Remark~7.2 there, with  the definition of \cite{AN} for $\Hdot^\beta_{v}$   that we adopt here.

  For  $\I=(0,\infty)$, the argument is mentioned in Step 4 of   the proof of  \cite[Theorem 5.2]{AIN}. It uses  the backward Kolmogorov operator  $\cK^{-}_{\beta, 0}$ with the same arguments as when $\I=\R$. Again, with our definition of Sobolev spaces when $\beta\ge d/2$, the hypothesis $f\in \L^2_{t,x}\Hdot^{\beta}_{\vphantom{t,x} v}$ suffices.    This would be the same on any interval $(T,\infty)$ by translating time. On intervals $(-\infty,T)$, it is obtained by  working with the forward Kolmogorov operator $\cK^+_{\beta,0}$ or by a change of variables  from the previous case. 
 
 \
  
 {\em Step 2:   $\ka>0$.} 
  
 The kinetic embeddings are given in  \cite[Theorem 6.1]{AIN} when $\ka=1$. For the general case when $\I$ is unbounded, one can follow the argument in \cite{AIN}. We give some details because it is instructive for the rest of the article to see how the parameter $\ka$ is handled. 
 
 When $\I=\R$, the argument combines  two ingredients. First, one proves  a uniqueness result in $\L^2_{t,x}\H^\beta_{v}$ for the corresponding equation: if $f\in\L^2_{t,x}\H^\beta_{v} $ solves  the distributional equality
 $$
 \pm (\partial_{t}+v\cdot\nabla_{x})f + (-\Delta_{v})^{\beta}f+\ka^{2}f=0$$
 then $f=0$. The argument is the same as  \cite[Lemma~3.1]{AIN}, and even easier because $\ka> 0$ and $f\in \L^2_{t,x,v}$. Next, one proves bounds for $\cK^+_{\beta,\ka}$ (and similarly  for $\cK^-_{\beta,\ka}$) on adapting the  situation where $\ka=0$ of  \cite[Section~2]{AIN}. The Kolmogorov kernel $K$ and the weight  $W$ there become
 \begin{equation}
  \label{eq:kolmogorovkernel}
K_{\beta,\ka}(t,s,\varphi,\xi)= \exp\bigg({-\int_{s}^t |\xi-\tau \varphi|^{2\beta}\, d\tau}- \ka^{2}(t-s)\bigg)
\end{equation}
and
\begin{align}
\label{eq:weight}
  W_{\beta,\ka}(t,\varphi,\xi)= \sup(|\xi-t \varphi|, |\varphi|^{\frac{1}{1+2\beta}}, \ka^{\frac 1 {\beta}}).
\end{align}
 With these new kernels and weights, the weighted operator estimates (bearing only on size estimates and Schur techniques)  can be proved identically with uniform estimates in $\ka$ (that is, $\ka$ does not appear in the operator bounds but only in the space norms),  so  that  we reach the conclusion of the statement for the dependency in $\ka$.  
 
 In the case of a half-line,  the same remarks as when $\ka=0$ apply. Further details are left to the reader.

  For the case of a bounded interval,  say $(0,T)$, we follow the localisation argument given in  \cite[Section 6.2]{AIN}. Assume  $f$ satisfies the hypotheses of the statement with $\ka^2=T^{-1}$. One decomposes 
  $f=\chi f+(1-\chi)f$ when $\chi:\R\to [0,1]$ is a smooth cut-off that is 1 on $(-\infty, T/4)$ and 0 on $(3T/4, \infty)$   with $\|\chi'\|_{\infty}\sim T^{-1}=\ka^2  $. Then, we can see that $\chi f$ and $(1-\chi)f$  satisfy the required assumptions to apply the conclusions on the intervals $(0,\infty)$ and $(-\infty,T)$ respectively. It remains to combine them.  
    \end{proof}

\begin{rem}
 Another way to follow the parameter $\ka$ is to reduce to $\ka=1$ by kinetic scaling when $\I$ is unbounded or when $\I$ is of size $\ka^{-2}$. 
 \end{rem}
 \begin{rem}
 It may look strange that $\ka$ appears in \eqref{eq:embedding} in front of $f$ but not of $S_{4}$. Let us explain this by considering the solution of the heat equation $(\partial_{t}-\Delta_{v}+\ka^2)f =\ka S_{4}$ given by 
 $f(t,v)=\int_{-\infty}^t e^{(t-s)(\Delta_{v}-\ka^2)}\ka S_{4}(s,v)\, \dd s$. Multiplying the equality by $\ka$ and using the contractivity of the heat semi-group on $\L^2_{v}$ we have $\ka\|f(t,\cdot)\|_{\L^2_{v}}\le \int_{-\infty}^t e^{-(t-s)\ka^2}\ka^2 \|S_{4}(s,\cdot)\|_{\L^2_{v}}\, \dd s$, so that Young's inequality for convolution in time yields
   $\ka\|f\|_{\L^2_{t,v}}\le \|S_{4}\|_{\L^2_{t,v}}$. \end{rem}

\section{Class of operators}\label{sec:thesetup}
We first define the class of operators $\cA$ with unbounded coefficients on $\Omega:=\R\times \,\R^{2d}$ as follows.  When we work on a strip $\I\times \,\R^{2d}$, where $\I$ is an open interval of $\R$, we implicitly assume that we have restricted the operator $\cA$ from full space. Note that at this stage $\beta>0$ and is not necessarily an integer. So  our assumptions apply to both differential and  integro-differential operators. 

\

\

\paragraph{\textbf{(B) Intermediate spaces $\Bdot$}} 
We assume that  $(\Bdot, \|\ \|_{\Bdot})$  is a reflexive normed space with the following five axioms:
\begin{enumerate} 
\item (qualitative containment)   $ \cD(\Omega) \subset \Bdot \subset \L^1_{t,x,\loc}\cD'_{v} 
$ where the notation $\cD'_{v}$ is for the space of distributions in the $v$-variable.
\item (canonical section) there exists a space $\Bdot^{\vphantom{'}}_{v}\subset \cD'_{v} $ such that for all $f\in \Bdot$ and almost every $(t,x)\in \R\times \, \R^d$, $f(t,x,\cdot)\in \Bdot_{v}$. 
\item (Stability for time multiplication) For all measurable subsets $\I\subset \R$, $f\in  \Bdot$ implies that  $f_{I}:=f\chi_{\I}$  belongs to $\Bdot$ and we have 
$\|f\chi_{\I}\|_{\Bdot}\le \|f\|_{\Bdot}$, where $\chi_{\I}$ is the indicator function of $\I$. 
\item (Extension by 0) 
Set $\Bdot_{\I}= \{f_{\mid\I}: f \in \Bdot\}$ equipped with the norm $\inf_{f\in \Bdot} \|f_{\mid\I}\|_{\Bdot}$ where $\I\subset \R$ is a measurable set and $f_{\mid \I}$ is the restriction to $\I$ of $f$. If  $g\in \Bdot_{\I}$ and $\tilde g$ is the extension   by zero  outside $\I\times \,\R^{2d}$ of $g$, then $\tilde g$ belongs to $\Bdot$ and $\|\tilde g\|_{\Bdot}= \|g\|_{\Bdot_{\I}}$.  

\item (embeddings)  There are  constants $0<c_{1} \le 1 \le c_{2}<\infty$ such that \begin{align}\label{eq:B}
(\Wbd)' \hookrightarrow& \  \Bdot \hookrightarrow \Fbd 
 \\
 \label{eq:B'}
(\Fbd)'  \hookrightarrow & \, \Bdot' \hookrightarrow \Wbd 
\end{align}
with 
\begin{align}
 \label{eq:embedB} c_{1}\|f\|_{\Fbd} \le&  \|f\|_{\Bdot}\le c_{2}\|f\|_{(\Wbd)'}\\
 \label{eq:embedB'} c_{2}^{-1} \|S\|_{\Wbd} \le & \|S\|_{\Bdot'} \le c_{1}^{-1} \|S\|_{(\Fbd)'} 
\end{align}
where $\Bdot'$ is the dual space of $\Bdot$ for the $\L^2_{t,x,v}$ inner product. 
\end{enumerate}

\

\paragraph{\textbf{(A) Operators $\cA$}}  
We set $\B=\Bd\cap \L^2_{t,x,v}$ with $\Bdot$ as in $\mathbf{(B)}$ and $\B^{\vphantom{2}}_{v}=\Bdot^{\vphantom{2}}_{v}\cap \L^2_{v}$. 
 We assume that for almost all $(t,x)\in \R\times \R^d$,  there are  sesquilinear forms  $$a_{t,x}\colon \B_{v} \times \B_{v} \to \IC$$ such that 
\begin{itemize}
\item $(t,x)\mapsto a_{t,x}(f,g)$ is measurable  for all $(f,g)\in \B_{v}\times \B_{v}\, , $ 
\item there are constants $\Lambda^{a}\ge 0$ and $c^{a}\ge 0$ such that
for all  $ f,g\in \B$ 
\begin{align} 
\label{eq:bounda}\iint_{\R\times \R^d }  |a_{t,x}(f(t,x,\cdot),g(t,x,\cdot))|\dt\dx \le  \Lambda^{a} (\|f\|_{\Bdot}+c^{a}\|f\|_{\L^2_{t,x,v}})( \|g\|_{\Bdot}+c^{a}\|g\|_{\L^2_{t,x,v}}).\end{align} 
\end{itemize}
We  set
\begin{equation}
  \label{e:a-defi}
  \forall f,g \in \B, \qquad a(f,g)=\iint_{\R\times \R^d} a_{t,x}(f(t,x, \cdot),g(t,x, \cdot))\, \dt\dx
\end{equation} 
and
 $\cA: \B \to \B'$ the continuous operator associated to $a$ by 
\begin{equation}
  \label{eq:cA} \angle{\cA f}{g}=  a(f,g), \quad f,g\in \B.\end{equation}
 We have 
\begin{equation}
\label{eq:upperboundcA}
\forall f,g \in \B, \quad  |\angle{\cA f}{g}|\le \Lambda^{a} (\|f\|_{\Bdot}+c^{a}\|f\|_{\L^2_{t,x,v}})( \|g\|_{\Bdot}+c^{a}\|g\|_{\L^2_{t,x,v}}).
\end{equation}

\begin{rem} Let us make some comments concerning $\mathbf{(B)}$.

\begin{enumerate}[label=\Alph*.]
\item Candidates for intermediate spaces are  mixed Lebesgue spaces, see the examples in Section~\ref{sec:examples}. 
 \item The conditions (ii), (iii), (iv) make sense thanks to the local integrability in time in (i).  Also the embeddings in (v) restrict to any subset $\I$ by the standard restriction/extension by zero retraction diagrams, and with the same constants $c_{1},c_{2}$.
 
 \item The space $\Fbd$ satisfies all items (i-v). It is a reflexive (Hilbert)  normed space  when $\beta<d/2$ and for $\beta\ge d/2$ it is only semi-normed. By abuse of terminology, we also say that $\mathbf{(B)}$ holds for $\Fbd$, but we shall avoid to use it when $\beta\ge d/2$.       
\item We use the dot notation $\dot{} $ to indicate that the norm is homogeneous with respect to the kinetic scaling $(t,x,v) \mapsto (\delta^{2\beta} t, \delta^{2\beta+1} x, \delta v)$. For $\Bd$ it is the same homogeneity as for $\Fbd$ and $(\Wbd)'$. For $\Bd'$ it is the one of $\Wbd$ and $(\Fbd)'$. Note that  the inclusions in \eqref{eq:B} and \eqref{eq:B'} are dense.
\item The embedding \eqref{eq:embedB} makes  $\Bdot$ a set of tempered distributions because the $\Fbd$ is constructed as such a subset, but the inclusion is not necessarily continuous. 
 
\end{enumerate}
 
 \end{rem}

 \begin{rem}
 The condition in  $\mathbf{(A)}$ is stable under adjunction of forms and operators. If $c^{a}=0$ then, we 
 have $\cA:\Bdot\to {\Bdot}'$. The condition \eqref{eq:bounda} is also stable by restriction to an interval $\I$ in time. 
\end{rem}

We record this result, seen from the  remarks above. 
 
\begin{lem}[The kinetic operator] Under  $\mathbf{(B)}$ and $\mathbf{(A)}$ then $\cA: \B \to \B'$ is bounded and 
 $(\partial_{t}+v\cdot \nabla_{x}) +\cA : \B \to \cD'(\Omega)$ is continuous. If, moreover, $c^{a}=0$, we may replace $\B$ by $\Bdot$. 
 \end{lem}

 To construct solutions we will often invert $(\partial_{t}+v\cdot \nabla_{x}) +\cA+\ka^2$ for large enough spectral parameter  $\ka\ge 0$.  In order to track  $\ka$  in the estimates we adopt the following notation.  
 
 For (solution) subspaces $\Bdot$ of $\Fbd$ as in $\mathbf{(B)}$, we define $\B_{\ka}$ as follows: \begin{itemize}
\item For $\ka=0$, we let $\B_{0}=\Bdot$ with norm $\|f\|_{\vphantom{\Bd}\B_{0}}=\|f\|_{\Bdot}$.
\item For $\ka>0$, we let $\B_{\ka}=\Bdot\cap \L^2_{t,x,v}$ with norm $\|f\|_{\vphantom{\Bd}\B_{\ka}}=\|f\|_{\Bdot}+\ka\|f\|_{\L^2_{t,x,v}}$.
\end{itemize}
If $\Bdot$ is a Hilbert space, we use instead the equivalent hilbertian norm $(\|f\|_{\Bdot}^2+\ka^2\|f\|_{\L^2_{t,x,v}}^2)^{1/2}$. 
  We define in this way $\Vbk$, $\Fbk$, $(\Wb)'_{\ka}$ and their associated norms.  
  
  For (source) subspaces $\Edot$ of  $\Wbd$ and $\ka\ge 0$, we define  $[\E]_{\ka}$ as follows:
   \begin{itemize}
\item For $\ka=0$, we let $[\E]_{0}=\Edot$ with norm $\|f\|_{\vphantom{\Edot}[\E]_{0}}=\|f\|_{\Edot}$.
\item For $\ka>0$, we let $[\E]_{\ka}$ be the space of $S_{1}+\ka S_{2}$ with $S_{1}\in \Edot$ and $S_{2}\in \L^2_{t,x,v}$ equipped with the norm $\|S\|_{\vphantom{\Edot}[\E]_{\ka}}=\inf (\|S_{1}\|_{\Edot}+\|S_{2}\|_{\L^2_{t,x,v}})$ (with modification if $\Edot$ is Hilbertian) where the infimum is taken over all such decompositions.  
\end{itemize}
   With this convention, 
 one can check that $([\Wb]_{\ka})'= {(\Wb)'}_{\ka}$, $(\B_{\ka})'=[\B']_{\ka}$ and $(\Fbk)'= [(\Fb)']_{\ka}$ with isometric norms when $\ka\ge 0$.  Recall that duality is understood with respect to duality extending the $\L^2_{t,x,v}$ inner product, explaining the absence of $\ka$ in the formula for the $[\E]_{\ka}$ norm. 
 
 Now that this is made precise, we omit the brackets $[]$ in the notation $[\E]_{\ka}$ to simplify the exposition.   We write $\Wbk$, $\B'_{\ka}$, etc. 
 
 If  the embeddings (v) in $\mathbf{(B)}$ hold, they extend to  the $\ka$ versions of the spaces with the same constants as we assumed $c_{1} \le 1 \le c_{2}$. 
 
 \

 \section{Weak solutions on full space: a priori properties}
 
Recall that $\Omega=\R\times \,\R^{2d}$.
 We are given $\Bdot, \cA$ as in  $\mathbf{(B)}$ and $\mathbf{(A)}$ and we are also given a spectral parameter $\ka\ge c^{a} $ with    $c^{a}$ in \eqref{eq:bounda}. 
 
 \begin{defn}[Weak solutions on $ \R\times \,\R^{2d}$]
\label{defn:weaksol}  
 Let $S \in \cD' (\Omega)$.  A  distribution $f \in \cD' (\Omega)$ is said to be a \emph{weak solution} to \begin{align}
\label{eq:CPR}
 (\partial_{t}+v\cdot\nabla_{x})f + \cA f+\ka^2 f= S,   
\end{align}
 if $f\in \B_{\ka}$ and for all $h\in \cD (\Omega
)$, 
\begin{align}
 \label{eq:CP2R}
  -\langle f, (\partial_{t}+v\cdot\nabla_{x})  h\rangle  + a(f,h) +\ka^2 \angle{f}{h} = \angle{S}{h}.
\end{align}
\end{defn}

\begin{rem}
 We interpret the first term as a sesquilinear pairing between a  distribution and a test function since $\Bka$ is a space of distributions.  This is an integral in $t,x$ of a distribution bracket in $v$. 
Whenever $\Bka$ is contained in $\L^1_{\loc}(\Omega)$, we have  
\[\langle f, (\partial_{t}+v\cdot\nabla_{x})  h\rangle=\int_{\Omega}  f \, (\partial_{t}+v\cdot\nabla_{x}) \overline h \, \dd(t,x,v).
\]
 Also  $\cD(\Omega) \subset \Bka$   since $ (\Wbk)' \subset \Bka$. The definition of $a(f,h)$ makes sense when $\ka=0$ as $\ka\ge c^{a}$ implies  $c^{a}=0$ and $\cA$ is defined on $\B_{0}$ in that case.
  
 \end{rem}

\begin{prop}[A priori regularity and estimates for weak solutions] \label{prop:aprioriestimates}
 Let  $f\in \B_{\ka}$ be a weak solution to \eqref{eq:CPR}  with $S\in \Wbk$. Then
\begin{enumerate}
\item $f$ belongs to $\Vbk$.
 \item $f$ satisfies  for all $-\infty\le \sigma\le \tau\le\infty$ the energy equality 
\begin{align}\label{eq:energyequality}
 \|f(\tau)\|_{\L^2_{x,v}}^2  + 2\Re \int_{\sigma}^\tau\int_{\R^d}(a_{t,x}(f,f)&+\ka^2\angle {f}{f})\, \dx\dt\\ \nonumber &=   \|f(\sigma)\|_{\L^2_{x,v}}^2+  2\Re \int_{\sigma}^\tau \angle {S(t)}{f(t)} \dt, 
\end{align}
where $\angle {S(t)}{f(t)}$ is to be interpreted as the term within parentheses of the right hand side in   \eqref{e:abs-cont}.

\item 
There are positive reals $\varepsilon_{0}=\varepsilon_{0}(d,\beta,\Lambda^{a}, c_{1},c_{2})$ and  $M=M(d,\beta,\Lambda^{a}, c_{1},c_{2})$, independent of $f$,  such that if $\ka \ge  c^{a}$ and if $\I$ is an unbounded interval, then we have
\begin{align}\label{eq:uniformA}
\| f_{\I} \|_{\L^2_{t,v}\Hdot_{\vphantom{t,v} x}^{\gb}}+ \sup_{t\in \overline{\I}} \|f(t)\|_{\L^2_{x,v}} 
\le M(\|f_{\I}\|_{\Bdot}
+  \ka \|f_{\I}\|_{\L^2_{\vphantom{t,x} t,x,v}} +\|S_{\I}\|_{\Wbk}).
 \end{align}
 In particular, with a different $M$ of the same nature, 
\begin{equation}\label{eq:estimateweaksolA}
\|f\|_{\Vbk}\le M( \|f\|_{\B^{\vphantom{\beta}}_{\ka}}+\|S\|_{\Wbk}).
\end{equation}
\end{enumerate}
        \end{prop}

\begin{rem}\label{rem:constants} 
\begin{itemize}
\item This is  a priori information about weak solutions and it does not use   coercivity assumptions.
\item The control of the $\L^2_{t,v}\Hdot_{\vphantom{t^a} x}^{\gb}$-norm  comes from the transfer of regularity and cannot be obtained from the energy inequality even when assuming coercivity of $\cA$.
\item  It is important to note that the constant $M$  does not depend on $c^{a}$, while the least possible  $\ka$ does. Moreover, $M$ is uniform in $\ka$.
\item  If  $c^{a}=0$,  then $\ka=0$ is possible and this gives a homogeneous estimate.

\end{itemize}

\end{rem}

\begin{proof} The items (i) and (ii) follow directly from Theorem~\ref{thm:optimalregularity}. 

We turn to the proof of (iii). We  prove \eqref{eq:uniformA} which, in particular, yields  the bound \eqref{eq:estimateweaksolA}. We assume $\I=\R$ and set $\L^2=\L^2_{t,x,v}$. As  $S\in \Wbk$,  it has a decomposition  $S=S_{1}+S_{2}+S_{3}+\ka S_{4}$ as in 
 Theorem \ref{thm:optimalregularity}, the norms of the $S_{i}$ being bounded independently of $\ka$ in the respective spaces. Writing  \begin{align*}
  (\partial_{t}+v\cdot\nabla_{x})f = S- \cA f-\ka^2 f =S- \cA f-\ka(\ka f),
\end{align*}
 Applying the inequality \eqref{eq:embedding}  and using $\|\ka(\ka f)\|_{\Wbk} \le   \|\ka f\|_{\L^2}$, we obtain 
\begin{align*}
 \| f \|_{\Gbd
 }+ \sup_{t\in {\R}} \|f(t)\|_{\L^2_{x,v}} 
&\lesssim_{d,\beta} \|f\|_{\Fbd
}
+  \ka \|f\|_{\L^2} +\|S-\cA f -\ka(\ka f)\|_{\Wbk},
\\
&\lesssim_{d,\beta} \|f\|_{\Fbd
}
+  \ka \|f\|_{\L^2} +\|S\|_{\Wbk}+\|\cA f \|_{\Wbk}.
\end{align*}
When  $\ka\ge c^{a}$, \eqref{eq:upperboundcA} yields for all $ f,g \in \B$,
\begin{align*}
  |\angle{\cA f}{g}|\le \Lambda^{a} (\|f\|_{\Bdot}+\ka \|f\|_{\L^2})\|g\|_{\Bka},\end{align*}
hence 
$$\|\cA f\|_{(\B_{\ka})'} \le  \Lambda^{a}  (\|f\|_{\Bdot}+\ka \|f\|_{\L^2})$$
and by \eqref{eq:embedB'}
$$\|\cA f\|_{\Wbk}\le  \Lambda^{a} c_{2} (\|f\|_{\Bdot}+\ka \|f\|_{\L^2}).$$
Thus, substituting this into the previous inequality and using \eqref{eq:embedB}, 
\begin{align*}
\| f \|_{\Gbd
 }+ \sup_{t\in {\R}} \|f(t)\|_{\L^2_{x,v}}
&\lesssim_{d,\beta} c_{1}^{-1}\|f\|_{\Bdot}+  \ka \|f\|_{\L^2} +\|S\|_{\Wbk}+\Lambda^{a} c_{2} (\|f\|_{\Bdot}+\ka \|f\|_{\L^2})
\\
&\lesssim_{d,\beta, \Lambda^{a}, c_{1},c_{2}} (\|f\|_{\Bdot}+  \ka \|f\|_{\L^2} +\|S\|_{\Wbk})
\end{align*}
as desired. 

The argument is rigorously the same if $\I$ is an unbounded interval by taking restrictions on $\I\times \,\R^{2d}$. 
\end{proof}
 
  \section{The duality scheme}
  \label{sec:duality}

Consider $\Bdot, \cA$ as in $\mathbf{(B)}$ and $\mathbf{(A)}$, and $\ka\ge c^{a}$. A priori regularity  suggests  that weak solutions constructed in 
$\B_{\ka}$ with appropriate sources belong to a better space.  We formalize this in the following definition. 
  
\begin{defn}[``Inversion'' of the kinetic operator]\label{def:inversion}    Let $\Z, \Y$  be  normed spaces. 
 We write $((\partial_{t}+v\cdot \nabla_{x}) +\cA+ \ka^2)^{-1}: \Z \to \Y$  if 
there exists a constant $C>0$ such that 
 \begin{enumerate}
 \item 
for each $S\in \Z$ there exists a  unique weak solution in $\B_{\ka}$  to 
 \eqref{eq:CPR}  and 
 \item    this solution belongs to $Y$ with  
$\|f\|_{\Y}\le C\|S\|_{\Z}$.
\end{enumerate}
We write   $\|((\partial_{t}+v\cdot \nabla_{x}) +\cA+ \ka^2)^{-1}\|_{\Z\to \Y} \le C$.

 There is the same definition for the backward equation.
 \end{defn}

Note that we are not saying in this definition that $(\partial_{t}+v\cdot \nabla_{x}) +\cA+ \ka^2$ is invertible from $\Y$ to $\Z$. We are just saying that there is the additional  bound  (ii) for the operator $((\partial_{t}+v\cdot \nabla_{x}) +\cA+ \ka^2)^{-1}: \Z \to \B_{\ka}$ defined by (i).
 In the definition,  $C$ can depend on $\ka$  but we will look for independence with respect to $\ka$ in a range.  
 
 The first result streamlines the duality scheme in \cite{AIN} and allows to reach optimal regularity of weak solutions with largest possible source space.  This is the first result in the article which does not always apply to the  $\ka=0$ case under our assumptions.

\begin{prop}[Duality scheme]\label{prop:duality} If $\ka>0$, the following assertions  are equivalent.
\begin{enumerate}
\item $\|((\partial_{t}+v\cdot \nabla_{x}) +\cA+ \ka^2)^{-1}\|_{\Bka'\to \Bka}  \le C_{1}$. 
\item $ \|(-(\partial_{t}+v\cdot \nabla_{x}) +\cA^*+ \ka^2)^{-1}\|_{\Bka'\to \Bka} \le C_{2}$. 
\item $\|((\partial_{t}+v\cdot \nabla_{x}) +\cA+ \ka^2)^{-1}\|_{\Wbk\to \Vbk} \le C_{3}$.
\item $ \|(-(\partial_{t}+v\cdot \nabla_{x}) +\cA^*+ \ka^2)^{-1}\|_{\Wbk\to \Vbk} \le C_{4}$.
 \end{enumerate}

 If $\ka=0$, then we have that both (iii) and (iv) are equivalent and imply (i) and (ii). 
 
 The relations between the constants depend only on $d,\beta,\Lambda^\cA, c_{1},c_{2}$.
\end{prop}

\begin{proof} In the case $\ka>0$, we show that (iii) $\Rightarrow$ (i) $\Rightarrow$ (iv) $\Rightarrow$ (ii) $\Rightarrow$ (iii). 
 If we assume (iii), then by embeddings in $\mathbf{(B)}$,
 \[\|((\partial_{t}+v\cdot \nabla_{x}) +\cA+ \ka^2)^{-1}\|_{\Bka'\to \Bka}   \le c_{2}^2C_{3}.\] 
 The proof of (iv) $\Rightarrow$ (ii) is the same.
 
 Next, we only prove (ii) $\Rightarrow$ (iii), the proof of (i) $\Rightarrow$ (iv) being identical.  So assume (ii). 
 We first establish  uniqueness for the  equation in (iii). The argument is valid for $\ka\ge 0$. So let $ f\in \Bka$ such that $((\partial_{t}+v\cdot \nabla_{x}) +\cA+ \ka^2) f=0$. First, by Proposition~\ref{prop:aprioriestimates}, we know that $ f\in \Vbk$ and $t\mapsto \|f(t)\|_{\L^2_{x,v}}^2$ is absolutely continuous on $\R$. By (ii), for $S\in \Bka'$ there exists $\tilde f \in \Bka$ such that 
$(-(\partial_{t}+v\cdot \nabla_{x}) +\cA^*+ \ka^2)\tilde f=S$. Again, $ \tilde f\in \Vbk$ and $t\mapsto \|\tilde  f(t)\|_{\L^2_{x,v}}^2$ is absolutely continuous on $\R$. 
Now, integrating the polarized version of \eqref{e:abs-cont}   and using the $\L^2_{x,v}$ limits 0 at $\pm\infty$,  with a slight abuse  in writing the dualities,  
\begin{align*}
0=[\langle  f(t), \tilde f(t)\rangle _{\L^2_{x,v}}]^{+\infty}_{-\infty}
=  \langle  -(\cA  +\ka^2) f, \tilde f \rangle
+ \langle  f, -S + (\cA^*+ \ka^2)\tilde f \rangle = -\langle  f, S\rangle.
\end{align*}
Hence, $ f=0$.

We next follow the duality scheme of \cite{AIN} to conclude for (iii) when $\ka>0$.  First,    we have  
$$(-(\partial_{t}+v\cdot \nabla_{x}) +\cA^*+ \ka^2)^{-1}:\B_{\ka}'\to \Vbk,
$$
with bound  $M(C_{2}+c_{2})$. 
Indeed, let $f=(-(\partial_{t}+v\cdot \nabla_{x}) +\cA^*+ \ka^2)^{-1}h$. Then, by \eqref{eq:estimateweaksolA},  (ii) and the embedding from $\mathbf{(B)}$, 
\begin{align*}
 \|f\|_{\Vbk}& \le M(\|f\|_{\Bka^{\vphantom{\beta}}}+\|h\|_{\Wbk})
 \\
 &  \le M (C_{2}\|h\|_{\Bka'^{\vphantom{\beta}}}+ \|h\|_{\Wbk})
\\
&\le M (C_{2} +c_{2} )\|h\|_{\Bka'^{\vphantom{\beta}}}.
\end{align*}
Introduce the continuous operator $\cT_{\ka}: \Wbk \to  \Bka$, since $\Bka$ is reflexive, defined by 
\[ \langle \cT_{\ka}S, h\rangle = \langle S, (-(\partial_{t}+v\cdot \nabla_{x}) +\cA^*+ \ka^2)^{-1}h\rangle,  \quad S\in \Wbk, h\in \Bka'.
\]
Its bound is not exceeding $M (C_{2} +c_{2} )$. We claim that if $S\in  \Bka'$ then $f=\cT_{\ka}S$ is a weak solution to  \eqref{eq:CPR}, that is $\cT_{\ka}S=((\partial_{t}+v\cdot \nabla_{x}) +\cA+ \ka^2)^{-1}S$ since uniqueness is granted. 
Indeed, $f\in \Bka$ by definition of $\cT_{\ka}$ and for $h\in \cD(\Omega)$,
\begin{align*}
-\langle f, (\partial_{t}+v\cdot\nabla_{x})  h\rangle  &+ a(f,h) +\ka^2 \angle{f}{h}\\
&= 
 -\langle \cT_{\ka}S, (\partial_{t}+v\cdot \nabla_{x})h\rangle + \langle \cT_{\ka}S, \cA^* h \rangle + \ka^2 \langle \cT_{\ka}S, h\rangle 
 \\
 &= \langle \cT_{\ka}S, -(\partial_{t}+v\cdot \nabla_{x})h+\cA^* h  +\ka^2 h \rangle
 \\
 &=\langle S, h\rangle.
\end{align*}
We have used  that $\cA^*+\ka^2:\Bka^{\vphantom{'}}\to \Bka'$ under our convention and also $\cD(\Omega)\subset \Bka'$ so the brackets from the second line   are all understood in the  $\Bka^{\vphantom{'}}, \Bka'$ duality and we can use the definition of $\cT_{\ka}$ and $h\in \Bka$ to get the last equality.

This claim holds for all  $S\in \Wbk$ by continuity of $\cT_{\ka}$  and  density of  $\Bka'$  in $\Wbk$ passing to the weak limit in the equation \eqref{eq:CP2R}. Thus we have obtained 
 $$((\partial_{t}+v\cdot \nabla_{x}) +\cA+ \ka^2)^{-1}:\Wbk\to \B_{\ka},
$$ with bound $M (C_{2} +c_{2} )$.   We use once more \eqref{eq:estimateweaksolA} to conclude 
$$((\partial_{t}+v\cdot \nabla_{x}) +\cA+ \ka^2)^{-1}:\Wbk\to \Vbk$$ 
with bound $M(M(C_{2}+c_{2})+c_{2})$. This proves (iii). 

Now we assume $\ka=0$. That (iii) implies (i) and (iv) implies (ii) is as before only using embeddings in $\mathbf{(B)}$. It remains to  establish that (iv) implies (iii) as the converse is the same. Assume therefore (iv). Then uniqueness for the forward equation in (iii) is the same argument as above.  For the existence in (iii), we also follow the argument above: we do not need to use \eqref{eq:estimateweaksolA} because $f\in \Vbd$ is part of the assumption. For the corresponding operator $\cT_{0}$, we   replace $\Bka'$ by $\Wbd$, $\Vbk$ by $\Vbd$, $\Wbk$ by $\Wbd$ and $\Bka$ by $(\Wbd)'$ and observe that $\cD(\Omega)\subset \L^1_{t}\L^2_{x,v}\subset  \Wbd$.
\end{proof}

\begin{rem} The equivalence between  (iii) and (iv) can be obtained with a direct proof when $\ka>0$ and one has  $C_{3}=C_{4}$.
 A slight variant of the argument gives a  direct proof that (i) and (ii) are equivalent with $C_{1}=C_{2}$ when $\ka>0$.    The main point when $\ka>0$ is the enlargement of the source space through this duality scheme while we already knew the weak solutions enjoy optimal regularity. 
\end{rem}

\begin{rem}
 If $\ka=0$, then the argument for existence in (ii) implies (iii) is valid provided $\cD(\Omega)\subset \Bdot_{0}'=\Bdot'$. Under the assumptions in $\mathbf{(B)}$, this is always granted if $\beta<d/2$ because it is the case of $(\Fbd)'$. But when $\beta\ge d/2$, this is not part of the assumption. If, for a specific choice of   $\Bdot$, it is true   then Proposition~\ref{prop:duality} holds for $\ka=0$. 
\end{rem}

\begin{prop}[The $\ka=0$ case as a limiting argument] \label{prop:limit} Any of the four items in Proposition~\ref{prop:duality} with uniform constant with respect to $\ka>0$ implies any of the four items with possibly different constants when $\ka=0$. 
\end{prop}

\begin{proof} As the four items  for $\ka>0$ are all equivalent, it suffices to prove that one of them for $\ka>0$ implies the same one for $\ka=0$. We do it for (i) as the three others have the same proof.

We begin with the existence. Let $S\in \Bdot'$. Then for all $\ka>0$,  $S\in \Bka'$ with 
$\|S\|_{\Bka'}\le \|S\|_{\Bdot'}$. By the assumption, there exists $f_{\ka}\in\Bka$ weak solution to $((\partial_{t}+v\cdot \nabla_{x}) +\cA+ \ka^2)f_{\ka}=S$ with $\|f_{\ka}\|_{\Bka}\le C_{1}\|S\|_{\Bka'}$.
By  \eqref{eq:estimateweaksolA}, $(\|f_{\ka}\|_{\Vbk})_{\ka>0}$ is a uniformly bounded family. Since $\|f_{\ka}\|_{\Vbd}\le \|f_{\ka}\|_{\Vbk}$, 
$(f_{\ka})_{\ka>0}$ is a uniformly bounded family in $\Vbd$. Hence it has a weak-star converging subsequence with $\ka_{n}\to 0$ in $(\Wbd)'$ to a weak solution $f$ to  $((\partial_{t}+v\cdot \nabla_{x}) +\cA)f=S$ (passing to the limit against a test function)  with  $\|f\|_{(\Wbd)'} \le C' \|S\|_{\Bdot'}$. Thus $\|f\|_{\Bdot} \le c_{2}C' \|S\|_{\Bdot'}$ by  \eqref{eq:embedB}.

We continue with uniqueness when $\ka=0$. First when (i) holds for $\ka>0$, we have (ii) for the same $\ka$ by Proposition~\ref{prop:duality}. Given $S\in \Bdot'$, the previous argument gives us existence of a weak solution  
$ \tilde f\in \Bdot$ to (ii) for $\ka=0$, namely $(-(\partial_{t}+v\cdot \nabla_{x}) +\cA^*)\tilde f=S$. If we assume that  $f\in \Bdot$ is a weak solution to $((\partial_{t}+v\cdot \nabla_{x}) +\cA) f=0$, then we can prove $\langle  f, S\rangle=0$ using the absolute continuity of $t\mapsto \langle  f(t), \tilde f(t)\rangle _{\L^2_{x,v}}$ and $\L^2_{x,v}$ limits at $\pm\infty$ as in  Proposition~\ref{prop:duality}.
\end{proof}

\begin{rem} $C_{1}$ is not preserved in this argument. The constants $C_{3}, C_{4}$ are preserved in the corresponding argument but not $C_{2}$ as well. However, if  $\cD(\Omega)\subset \Bdot'$ then one can proceed directly by weak star convergence  in $\Bdot$ without resorting to \eqref{eq:estimateweaksolA} and this  preserves $C_{1}$. The same remark holds for  $C_{2}$.
 \end{rem}

Let us state a corollary giving a form of maximal regularity.  For $\ka\ge 0$, define \begin{align*}
  \cB_{\ka}&=\{f\in \B_{\ka}; (\partial_{t}+v\cdot \nabla_{x})f\in \B_{\ka}'\}, \\
   \widetilde\cB_{\ka}&=\{f\in \B_{\ka}; (\partial_{t}+v\cdot \nabla_{x})f\in \Wbk\}, 
\end{align*}
 with norms
 \begin{align*}
 \|f\|_{ \cB_{\ka}}&=\|f\|_{\B_{\ka}}+ \|(\partial_{t}+v\cdot \nabla_{x})f\|_{\B_{\ka}'},
 \\
  \|f\|_{\widetilde \cB_{\ka}}&=\|f\|_{\B_{\ka}}+ \|(\partial_{t}+v\cdot \nabla_{x})f\|_{\Wbk}.
  \end{align*}
  Note that $\cB_{\ka}^{\vphantom{\beta}}\subset \widetilde \cB_{\ka}^{\vphantom{\beta}}$. 

\begin{cor}[maximal regularity]
\label{cor:isomorphism} We have the following maximal regularity estimates
\begin{enumerate}
\item Assume $((\partial_{t}+v\cdot \nabla_{x}) +\cA+ \ka^2)^{-1} : {\Bka'\to \Bka^{\vphantom{'}}}  $.  Then
 $(\partial_{t}+v\cdot \nabla_{x}) +\cA+ \ka^2: \cB_{\ka}^{\vphantom{'}} \to \B_{\ka}'$  is an isomorphism   with
\begin{align*}
 \|((\partial_{t}+v\cdot \nabla_{x}) +\cA+ \ka^2)f\|_{\B_{\ka}'} &\sim \|f\|_{\Bka} +  \|(\partial_{t}+v\cdot \nabla_{x})f\|_{\B_{\ka}'}.
\end{align*}
The space  $\cB_{\ka}$ is a Banach space and we have the embedding $\cB_{\ka}^{\vphantom{\beta}}\subset  \Vbk$ with implicit constant independent of $\ka$. 
\item Assume $((\partial_{t}+v\cdot \nabla_{x}) +\cA+ \ka^2)^{-1} : {\Wbk\to \Vbk}  $. Then
 $(\partial_{t}+v\cdot \nabla_{x}) +\cA+ \ka^2: \cB_{\ka}^{\vphantom{'}} \to \Wbk$  is an isomorphism  with   
\begin{align*}
 \|((\partial_{t}+v\cdot \nabla_{x}) +\cA+ \ka^2)f\|_{\Wbk} &\sim \|f\|_{\Bka^{\vphantom{\beta}}} +  \|(\partial_{t}+v\cdot \nabla_{x})f\|_{\Wbk}.
\end{align*}
The space $\widetilde \cB_{\ka}$ is a Banach space and we have the embedding $ \widetilde \cB_{\ka}^{\vphantom{\beta}}\subset \Vbk$ with implicit constant independent of $\ka$. 
\end{enumerate}
If the implicit constants do not depend on $\ka$  in the hypothesis, the same happens for the conclusions.
 \end{cor}
 
\begin{rem}
 It follows from the observations above that (i) and (ii) hold simultaneously when $\ka>0$ or when $\ka=0$ and $\cD(\Omega)\subset \Bdot'$. 
\end{rem}

\begin{proof} The proofs of (i) and (ii) are the same. By definitions of $\cB_{\ka}$ and $\widetilde \cB_{\ka}$,  the operators in the statement are bounded.  Injectivity comes from assumption. Surjectivity also follows from assumption since $((\partial_{t}+v\cdot \nabla_{x}) +\cA+ \ka^2)^{-1}$ maps into weak solutions to the equation. Thus bijectivity is proved and the continuity bounds follow immediately. The embedding is also clear from \eqref{eq:estimateweaksolA}.
\end{proof}

\begin{rem}\label{rem:homogeneous} The spaces $\cB_{\ka}$ and $\widetilde \cB_{\ka}$ are independent of $\cA$. The conclusions apply whenever there exists one  operator given by $\mathbf{(A)}$ with the assumptions. In the case   
   $\ka=0$, the hypotheses and conclusions  in each statement are  homogeneous.   \end{rem}

\section{The  Cauchy problem} 

Consider  again $\Bdot, \cA$ as in $\mathbf{(B)}$ and $\mathbf{(A)}$, and $\ka\ge c^{a}$.

Causality is the fact that a source being 0 before some time $t$ implies the weak solution vanishes before $t$. We say that $(\partial_{t}+v\cdot \nabla_{x}) +\cA+ \ka^2$ is causal if whenever $S\in \Wbk$ and  $\I$ is an open interval that is a neighbourhood of $-\infty$ with  $S=0$ on $\I\times \,\R^{2d}$, then the corresponding solution $f=((\partial_{t}+v\cdot \nabla_{x}) +\cA+ \ka^2)^{-1}S$ also vanishes on $\I\times \,\R^{2d}$. For the backward operator, change $-\infty$ to $+\infty$ and use anti-causal terminology.

\subsection{$\L^2$ valued measure sources} 
 
In preparation for the Cauchy  problem, we show the following principle: the   solvability estimate $$((\partial_{t}+v\cdot \nabla_{x}) +\cA+ \ka^2)^{-1}:\L^1(\R;\L^2_{x,v}) \to \Vbk$$
extends automatically to
$$((\partial_{t}+v\cdot \nabla_{x}) +\cA+ \ka^2)^{-1}:\cM(\R;\L^2_{x,v}) \to (\Wbk)'$$
where, with $\H=\L^2_{x,v}$, $\cM(\R;\H)$ denotes the space of bounded $\H$-valued measures on $\R$ defined as the topological dual space of $\C_0(\R; \H)$ with respect to the sup-norm. 

 We denote by $\llangle \mu, \varphi \rrangle_{\cM,\C_0}$ the sesquilinear duality bracket. We equip the space $\cM(\R;\H)$ with the induced norm, that is 
\begin{equation*}
    \left \|\mu   \right \|_{\mathcal{M}(\R;\H)}:= \sup_{\varphi \in \C_0(\R;\H)\setminus\left \{ 0 \right \}} \frac{\left | \llangle \mu, \varphi \rrangle_{\mathcal{M},\C_0} \right |}{\ \ \ \ \left \|\varphi  \right \|_{\L^\infty(\R;\H)}}.
\end{equation*}
It is a Banach space containing a subspace isometric to $\L^1(\R;\H)$.

For $\mu \in \cM(\R;\H)$, $\mathrm{supp}(\mu)$ is the complement of the largest open set of $\R$ on which $\mu$ is equal to $0$. More precisely, we say that $\mu$ equals $0$ on an open set $O \subset \R$ if for all $\varphi \in \C_0(\R;\H)$ with support contained in  $O$, $\llangle \mu, \varphi \rrangle_{\mathcal{M},C_0}=0$. Let $\mathcal{N}$ denote the set of all such open sets. We have 
\begin{align*}
    \mathrm{supp}(\mu ):= \R\setminus \bigcup_{O \in \mathcal{N}} O.
\end{align*}
 Dirac measures are important examples.  For any $s \in \R$ and $\psi \in \H$, we denote by $\delta_s \otimes \psi$ the Dirac measure at $s$ carried by $\psi$ which is defined by 
\begin{equation*}
    \llangle \delta_s \otimes \psi,  \varphi \rrangle_{\cM,\C_0} =  \langle \psi,  \varphi(s) \rangle_{\H}, \ \varphi \in \C_0(\R;\H).
\end{equation*}

 We first state the classical lemma below for later use.
\begin{lem}\label{Radon}
    Let $\mu \in \cM(\R;\H).$ Then there exists a sequence $(S_n)$ of elements in $\L^1(\R;\H)$ such that 
    $S_n \rightharpoonup \mu$  (weak-$\star$  convergence)  and  $\sup_{n} \left \| S_n  \right \|_{\L^1(\R;\H)} \leq \left \| \mu  \right \|_{\cM(\R;\H)}$.
   \end{lem}
\begin{proof}
We obtain the sequence $(S_n)$ by convolving $\mu$ in the $t$-variable with a scalar mollifying sequence $(\theta_n)$ and we easily check  that   the required properties hold.
\end{proof}

\begin{prop}[Existence and uniqueness for bounded measure source] \label{CasRadon} 
 Assume 
$$((\partial_{t}+v\cdot \nabla_{x}) +\cA+ \ka^2)^{-1}:\L^1(\R;\H) \to \Vbk.$$ Then
     $$((\partial_{t}+v\cdot \nabla_{x}) +\cA+ \ka^2)^{-1}:\cM(\R;\H) \to (\Wbk)'$$  with the same bound. 
 Let $\mu \in \mathcal{M}(\R;\H)$ and $f=((\partial_{t}+v\cdot \nabla_{x}) +\cA+ \ka^2)^{-1}\mu$.     If $\I \subset \R\setminus \mathrm{supp}(\mu)$ is an unbounded open interval, then the restriction of $f$ to $\I\times \,\R^{2d}$  extends to a function  $g\in \C_0(\overline{\I},\H)$ with $t \mapsto \left \| g(t) \right \|_{\H}^2$ is absolutely continuous on $\overline{\I}$. Moreover, if $(\partial_{t}+v\cdot \nabla_{x}) +\cA+ \ka^2$ is causal, then $f=0$ on ${\I}\times \,\R^{2d}$ for any $\I$ neighborhood of $-\infty$. 
\end{prop}

\begin{proof}
    Uniqueness in $\B_{\ka}$ follows from the hypothesis using (i) in definition~\ref{def:inversion}. To prove the existence, we pick $\left ( S_n \right )
        $ as in  lemma~\ref{Radon}. 
By assumption, $f_{n}=((\partial_{t}+v\cdot \nabla_{x}) +\cA+ \ka^2)^{-1}S_{n} \in \Vbk\subset (\Wbk)'$ with uniform bound with respect to $n$. Using the Banach-Alaoglu theorem in $(\Wbk)'$, and passing to the limit in the equation, there exists $f \in   (\Wbk)'$ weak solution to $(\partial_{t}+v\cdot \nabla_{x})f +\cA f+ \ka^2 f =\mu$ with 
\[ 
\|f\|_{(\Wbk)'} \le \liminf  \|f_{n}\|_{\Vbk} \lesssim C \liminf  \|S_{n}\|_{\L^1(\R;\H)} \le C\|\mu\|_{\cM(\R;\H)}.
\]

Now, if $\I \subset \R$ is an unbounded open interval such that $ \I \cap \, \mathrm{supp}(\mu)= \varnothing$, then by Theorem  \ref{thm:optimalregularity}, the restriction of $f$ to $\I\times \,\R^{2d}$ has an extension $g \in \C^{}_0(\overline{\I};\H)$ and $t \mapsto \left \| g(t) \right \|_{\H}^2$ is absolutely continuous on $\overline{\I}$. 

If $\I=(-\infty,T)$, we can see that the approximations $S_{n}$ vanish on any $(-\infty,T')\times \,\R^{2d}$ for any $T'<T$ and $n$ large enough. As $ \L^1(\R;\H)\subset \Wbk$,  we have $f_{n}=0$ on $(-\infty,T')\times \,\R^{2d}$  by the causality assumption. By taking subsequences, we conclude the same thing for the weak-star limit $f$, and  letting $T'\to T$ using $f\in \C^{}_0((-\infty,T];\H)$ concludes the argument.  
 \end{proof}

The next corollary is crucial for formulating the initial value problem and  constructing a fundamental solution.

\begin{cor}[Existence and uniqueness for Dirac measure source]\label{cor:initialvalue} Assume $((\partial_{t}+v\cdot \nabla_{x}) +\cA+ \ka^2)^{-1}:\L^1(\R;\H) \to \Vbk$ and  $(\partial_{t}+v\cdot \nabla_{x}) +\cA+  \ka^2$  is causal.
Let $s \in \R$ and $\psi \in \H$. There exists a unique $f \in \B_{\ka}$ solution to
   $(\partial_{t}+v\cdot \nabla_{x})f +\cA f+ \ka^2 f = \delta_s \otimes \psi$ in $\mathcal{D}'(\Omega).$
Moreover, $f\in (\Wbk)'$ with  
$$\|f\|_{(\Wbk)'}\le \|((\partial_{t}+v\cdot \nabla_{x}) +\cA+ \ka^2)^{-1}\|_{\L^1(\R;\H) \to \Vbk}\ \|\psi\|_{\H}\ ,$$
  $f\in \C_0(\R\setminus  \{ s  \}; \H)$, $f$ equals  $0$ on $(-\infty,s)$ and $\lim_{t \to s^+} f(t)=\psi$ in $\H$, and $t\mapsto \|f(t)\|_{\H}^2$ is absolutely continuous on $(s,\infty)$. 
\end{cor}

\begin{proof} Applying the previous result, we have existence in $(\Wbk)'$ and uniqueness in $\Bka$ of a weak solution  $f$ with   $f \in \C_0(\R\setminus  \{ s  \};\H)$, $f$ equals  $0$ on $(-\infty,s)$ and has a limit $f(s^{+})$ when  $t \to s^+$ and $t\mapsto \|f(t)\|_{\H}^2$ is absolutely continuous on $(s,\infty)$.

It remains to show $f(s^{+})=\psi$. Let $\Tilde{\psi}\in \cD(\R^{2d})$ and $\theta\in \mathcal{D}(\R)$ with $\theta(s)=1$ and set $\varphi=\theta\otimes \Tilde{\psi}$. Using the absolute continuity of $t\mapsto \langle f(t), \varphi(t)\rangle_{\H}$ on both $[s,\infty)$ and $(-\infty,s]$ following from  a polarisation argument, we have
\begin{align*}
     \langle f(s^+), \Tilde{\psi}\rangle_{\H} &= \int_s^\infty\int_{\R^d}  a_{t,x}(f,\varphi) \  \dt  \dx - \int_s^\infty - \ka^2\langle f, \varphi\rangle_{\H}  +
 \langle f, (\partial_t+v\cdot\nabla_{x})\varphi\rangle_{\H} \ \dt,    \\
    - \langle f(s^-), \Tilde{\psi}\rangle_{\H} &=  \int_{-\infty}^s\int_{\R^d}  a_{t,x}(f,\varphi) \  \dt  \dx - \int_{-\infty}^s - \ka^2\langle f, \varphi\rangle_{\H}  +
 \langle f, (\partial_t+v\cdot\nabla_{x})\varphi\rangle_{\H} \ \dt.   
    \end{align*}
    By the equation for $f$ on $\Omega$, the sum of the right hand sides is exactly  $$-\llangle \delta_s \otimes \psi,  \varphi \rrangle_{\cM,\C_0}= -
    \langle \psi, \Tilde{\psi}\rangle_{\H}   .$$ 
    By density of $\cD(\R^{2d})$ in $\H$, we conclude that $
 f(s^+)-f(s^-)= \psi$. Since we know that $f(s^-)=0$ by causality, we have proved $f(s^+)=\psi$. 
\end{proof}

\subsection{Cauchy problems}

We introduce the following condition. 

\

\paragraph{\textbf{Assumption (WPC): Well-posedness and causality}} 
We are given $\Bdot, \cA$ as in $\mathbf{(B)}$ and $\mathbf{(A)}$. We assume there exists  $\ka\ge c^{a}$ such that  $((\partial_{t}+v\cdot \nabla_{x}) +\cA+ \ka^2)^{-1}:\Wbk \to \Vbk$ and  $(\partial_{t}+v\cdot \nabla_{x}) +\cA+ \ka^2$ is causal.  

\
\begin{rem}
 The  condition $((\partial_{t}+v\cdot \nabla_{x}) +\cA+ \ka^2)^{-1}:\Bka'\to \Bka^{\vphantom{'}}$ is equivalent when $\ka>0$ by Proposition~\ref{prop:duality} but could be weaker when $\ka=0$.  In our examples, we prove them with $\ka>0$ with uniform bounds and pass to the limit. 
\end{rem}

Recall that the hypotheses in $\mathbf{(A)}$ restrict to  any open strip $\I\times \,\R^{2d}$ with the same constants (take $f,g$ supported in $\I$). 
When $\E$ is a space of distributions on $\R\times \,\R^{2d}$ and $\I$ is an open interval, we use the notation $\E(\I)$ to denote its restriction to $\I\times \,\R^{2d}$.

 \begin{defn}[Weak solutions to the kinetic Cauchy problem] \label{defn:weaksol-cauchy}  Let $0<T\le \infty$.
 Let $S \in \cD' ((0,T)\times \,\R^{2d})$  and $\psi \in \cD'(\R^{2d})$. Assume that the  distribution $S$ extends to a continuous linear functional on $\cD ([0,T)\times \,\R^{2d})$. A  distribution $f \in \cD' ((0,T)\times \,\R^{2d})$ is said to be a \emph{weak solution} to \begin{align}
 \label{eq:CP}
 \begin{cases}
 (\partial_{t}+v\cdot\nabla_{x})f + \cA f+\ka^2 f= S,  & \\
  f(0)=\psi,  &
\end{cases}
\end{align}
 if $f\in \Bka((0,T))$ and for all $h\in \cD ([0,T)\times \,\R^{2d})$, 
\begin{align}
 \label{eq:CP2}
  -\langle f ,  (\partial_{t}+v\cdot\nabla_{x})  h \rangle  + a(f,h) +\ka^2\angle f h  = \angle{S}{h} + \angle{\psi}{h(0)}.
\end{align}

\end{defn}  

With this notion, we have the following results.

\begin{thm}[Kinetic Cauchy problem on a half-space]\label{thm:CPhalfspace} Assume $\mathbf{(WPC)}$. Let $S\in \Wbk((0,\infty))$  and $\psi\in \L^2_{x,v}$. There exists a unique   weak solution to  \eqref{eq:CP} in $\Bka((0,\infty))$. Moreover,   $f\in \Vbk((0,\infty))$ with $f\in \C^{}_{0}([0,\infty)\, ;\, \L^2_{x,v}) $ and $ f(0,\cdot)=\psi$, and 
$$
\|f\|_{\Vbk((0,\infty))} \le\|((\partial_{t}+v\cdot \nabla_{x}) +\cA+ \ka^2)^{-1}\|_{\Wbk \to \Vbk}\ (\|S\|_{\Wbk((0,\infty))}+ \|\psi\|_{\L^2_{x,v}}).
$$
 \end{thm}
 
 \begin{thm}[Kinetic Cauchy problem on a strip]\label{thm:CPstrip} Assume $\mathbf{(WPC)}$. Let $0<T<\infty$,  $S\in \Wb((0,T))$ and $\psi\in \L^2_{x,v}$. There exists a unique  weak solution to  \begin{align*}
 \begin{cases}
 (\partial_{t}+v\cdot\nabla_{x})f + \cA f = S,  & \\
  f(0)=\psi,  &
\end{cases}
\end{align*}
 in  $\B((0,T))$. Moreover,  $f\in \Vb((0,T))$ with $f\in \C([0,T]\, ;\, \L^2_{x,v}) $ and  $f(0,\cdot)=\psi$,  
and
$$
\|f\|_{\Vbk((0,T))} \le \|((\partial_{t}+v\cdot \nabla_{x}) +\cA+ \ka^2)^{-1}\|_{\Wbk \to \Vbk}\ e^{T\ka^2} (\|S\|_{\Wbk((0,T))}+ \|\psi\|_{\L^2_{x,v}}).
$$
 \end{thm}
 
 In the last statement,  we have followed the spectral parameter $\ka$ to get the bound.  
 
\begin{proof}[Proof of Theorem~\ref{thm:CPhalfspace}]
For uniqueness we proceed as follows. First, by  Theorem~\ref{thm:optimalregularity}, we know that solutions belong to $\Vbk((0,\infty))$ and $\C^{}_{0}([0,\infty)\, ;\, \L^2_{x,v})$. Assume that $S=0$ and $\psi=0$. Then the extension $\tilde f$ by 0 of $f$ on $(-\infty,0)\times \,\R^{2d}$ belongs to $\Vbk(\R)=\Vbk$ since $\tilde f(0)=\psi=0$. Computing in $\cD'(\R\times \,\R^{2d})$ and using the equation for $f$, we find 
$((\partial_{t}+v\cdot\nabla_{x}) + \cA +\ka^2)\tilde f =\delta_{0}\otimes (f(0^+)- \tilde f(0^-))=0$. We conclude that $\tilde f=0$ by uniqueness in $\R\times \,\R^{2d}$. 

For existence, we pick $S\in \Wbk((0,\infty))$ and $\psi\in \L^2_{x,v}$. We let $S_{0}$ be the  zero extension  to  $\R\times \,\R^{2d}$ of $S$, and $f_{0}= ((\partial_{t}+v\cdot\nabla_{x}) + \cA +\ka^2)^{-1}S_{0}$ restricted to $(0,\infty)\times \,\R^{2d}$. We let $f_{1}$ be the restriction to $(0,\infty)\times \,\R^{2d}$ of the solution obtained in Corollary~\ref{cor:initialvalue} with data $\delta_{0}\otimes \psi$. We set $f=f_{0}+f_{1}$. The verification of \eqref{eq:CP2} follows on calculations using the absolute continuity of $t\mapsto \angle{f(t)}{h(t)}$ on $(0,\infty)$. 

The bound follows from combining  the preceding results. 
\end{proof}

\begin{proof}[Proof of Theorem~\ref{thm:CPstrip}]
First, by the usual exponential trick $f\to e^{-t\ka^2}f$,  the problem is equivalent to the one with $\cA+\ka^2$  with same initial data $\psi$ and $e^{-t\ka^2}S$ instead of $S$. 

 For uniqueness we proceed also by extension. First, by  Theorem~\ref{thm:optimalregularity}, we know that solutions belong to $\Vbk((0,T))$ and $\C^{}([0,T]\, ;\, \L^2_{x,v})$. Assume that $S=0$ and $\psi=0$. Consider the extension $\tilde f$ by 0 of $f$ on $(-\infty,0)\times \,\R^{2d}$, and by the solution $f^\sharp$ on $(T,\infty)\times \,\R^{2d}$ obtained in Corollary~\ref{cor:initialvalue} with data $\delta_{T}\otimes f(T)$. Then $\tilde f$  belongs to $\Vbk(\R)=\Vbk$ since $\tilde f(0)=\psi=0$ and  
 $f(T)=f^\sharp(T^+)$. Computing in $\cD'(\R\times \,\R^{2d})$ and using the equations for $f$ and $f^\sharp$, we find 
$$((\partial_{t}+v\cdot\nabla_{x}) + \cA +\ka^2)\tilde f =\delta_{0}\otimes ( f(0^+)- \tilde f(0^-))+ \delta_{T}\otimes (f(T^-)-f^\sharp(T^+))=0.$$ We conclude that $\tilde f=0$ by uniqueness in $\R\times \,\R^{2d}$.

For the existence, we simply restrict the solution obtained in the proof of Theorem~\ref{thm:CPhalfspace}. 

The bound follows by multiplying by $e^{\ka^2 t}$ at the end. 
\end{proof}

\section{The fundamental solution}

In this generality, there is an appropriate notion of a fundamental solution at the operator level.
Existence and uniqueness of a fundamental solution to $(\partial_{t}+v\cdot\nabla_{x}) + \cA$ was proved in \cite{AIN2}, Section 3, in the specific case where $\Bdot=\Fbd$ and $\cA$ has bounded coefficients. One can observe that this proof proceeds from an automatic functional analysis principle which generalizes to   assumption  $\mathbf{(WPC)}$. 

It is simpler to work on $\I=(0,T)$ with  $0<T<\infty$. We could also consider infinite intervals at the expense of technical adjustments. Alternatively, one can restrict to any finite intervals when working on infinite intervals. 

\begin{thm}[Existence and uniqueness of a fundamental solution]\label{thm:fso}
 Assume that $\mathbf{(WPC)}$ holds. Then $(\partial_{t}+v\cdot\nabla_{x}) + \cA$ possesses a unique fundamental solution on $[0,T]\times \, \R^{2d}$.  
\end{thm}

Fundamental solutions provide a representation of weak solutions.  For the sake of completeness, we recall the definition of a fundamental solution, summarize the functional analysis principle behind this theorem and review representation. Details are found in \cite{AIN2}.

\begin{defn}[Fundamental solution]
 \label{def:FSOF}
 The fundamental solution to $(\partial_{t}+v\cdot\nabla_{x}) + \cA$  is a family of 
 operators $\Gamma=(\Gamma(t,s))_{ (s, t)\in [0,T]^2}$, called fundamental solution operators, with the following properties.
 \begin{enumerate}
 \item They are uniformly bounded operators on $\L^2_{x,v}$.
  \item $\Gamma(t,s)=0$  for almost every $(s,t)$ with $s>t$.
 \item  For any $\psi,\tpsi \in \cD(\R^{2d})$,  $(s,t)\mapsto \angle {\Gamma(t,s)\psi}{\tpsi}$ is a Borel measurable function.  
 \item For any tensor product of test functions $S=\phi\otimes \psi$, the unique weak solution $f\in \B((0,T))$ to  the kinetic Cauchy problem 
   $(\partial_{t}+v\cdot\nabla_{x})f + \cA f=S$ on $(0,T)$  with $f(0)=0$, satisfies for all $\tpsi\in \cD(\R^{2d})$ and almost every $t\in (0,T]$ 
   \begin{align}
   \label{eq:representationfs}
  \angle{f(t)}{\tpsi}=\int_{(0,t)} \phi(s)\angle {\Gamma(t,s)\psi}{\tpsi}\, \dd s.
 \end{align}
  
  \end{enumerate}

  One defines  the fundamental solution $\widetilde \Gamma=(\widetilde \Gamma(s,t))_{ (s,t) \in [0,T]^2}$  to the backward operator   $-(\partial_{s}+v\cdot \nabla_{x})+\cA^*$ analogously  and (ii) is replaced by  $\widetilde \Gamma(s,t)=0$ if $s>t$.
  \end{defn}
  
  Note that from (i), (ii) and (iii),  $(s,t)\mapsto \angle {\Gamma(t,s)\psi}{\tpsi}$ is  integrable, that is,  
  \[\int_{(s,t)\in [0,T]\times [0,T]} |\angle {\Gamma(t,s)\psi}{\tpsi}|\, \dd s\dt <\infty,
   \]
   so that the integral in (iv) makes sense for almost every $t$. 
   
   \

 A fundamental solution, if it exists, must be unique up to almost everywhere equality. It is shown that the Green operator family is the  canonical representative. This family is given by the maps $G(t,s)$, where  $G(t,s)\psi$ is the solution at time $t$ in $\L^2_{x,v}$ of the initial Cauchy problem \eqref{eq:CP}   when $\ka=0$ with no source and data $\psi\in \L^2_{x,v}$ at time $s$ on $(s,T)\times \,\R^{2d}$ if $0\le s\le t\le T$. This is why we needed to treat the initial Cauchy problem. Thus we identify $\Gamma(t,s)$ with the Green operator $G(t,s)$. Hence, it  inherits all its properties  among which
 \begin{enumerate}
 \item   everywhere definition  and  boundedness in $\L^2_{x,v}$,
\item  continuity in the space of bounded operators on  $\L^2_{x,v}$ in the variables $t$ and $s$ separately (strong topology for $t$ and weak topology for $s$),
 \item  the adjoint relation $\Gamma(t,s)^*=\widetilde \Gamma(s,t)$ for all $s,t$,
 \item  the Chapman-Kolmogorov relation 
 $\Gamma(t,s)=\Gamma(t,r)\Gamma(r,s)$ for all $s\le r\le t$. 
\end{enumerate}
 Additional properties  come from the $\Vb$ estimates of weak solutions.  
 
 \
 
 We state, without proof, the representation theorem.
  
  \begin{thm}[Representation of weak solutions to the kinetic Cauchy problem] \label{thm:representation0T}
 Assume that $\mathbf{(WPC)}$ holds. Let $0<T<\infty$. For $S\in \Wb((0,T))$ and $\psi\in \L^2_{x,v}$, 
  the  unique  weak solution $f\in \B((0,T))$ given in Theorem~\ref{thm:CPstrip}  is represented  for any $t\in [0,T]$ in $\L^2_{x,v}$  by
  \begin{align}
 \label{eq:representationweaksolCP0T}
 f(t)=\Gamma(t,0)\psi+ \int_{0}^t \Gamma(t,s) S(s)\ds
 \end{align}
 where the integral is understood in the following sense.

 If $S\in \L^1_{t}\L^2_{x,v}$, then the integral is defined in the strong (Bochner) sense. 

 If $S= D_{v}^\beta\, F$ with $F \in  \L^2_{t,x,v}$, the integral is defined weakly in $\L^2_{x,v}$ by $($with $\L^2_{x,v}$ inner products$)$
 \begin{align}
 \label{eq:representationweaksol1}
 \angle{f(t)}{\tpsi}=\int_{0}^t \angle{F(s)}{D_{v}^\beta\, \widetilde\Gamma(s,t)\tpsi}\ds, \quad \tpsi \in \L^2_{x,v}.
 \end{align}

 If $S= D_{\vphantom{t}x}^{\gb} F$ with $F \in  \L^2_{t,x,v}$, the integral is defined weakly in $\L^2_{x,v}$ by $($with $\L^2_{x,v}$ inner products$)$
 \begin{align}
 \label{eq:representationweaksol2}
 \angle{f(t)}{\tpsi}=\int_{0}^t \angle{F(s)}{D_{x}^{\gb}\, \widetilde\Gamma(s,t)\tpsi}\ds, \quad \tpsi \in \L^2_{x,v}.
 \end{align}

 The analogous  representation  holds for  weak solutions to the backward kinetic Cau\-chy problem  for 
 $-(\partial_{t}+v\cdot \nabla_{x})+\cA^*$ with $\widetilde \Gamma(s,t)$. 
 \end{thm}

 Note that if we replace $\cA$ by $\cA+\ka^2$ and let $\ka$ vary, then uniqueness implies a relation between the different fundamental solutions, namely $$\Gamma_{\cA+\ka^2 }(t,s)=e^{-\ka^2(t-s)}\Gamma_{\cA}(t,s)$$ with obvious notation.

There are of course   corresponding  statements for weak solutions to  Cauchy problems on a half-space  for \eqref{eq:CP} and for  weak solutions to \eqref{eq:CPR} on full space for $(\partial_{t}+v\cdot\nabla_{x}) + \cA +\ka^2$ provided we assume $\mathbf{(WPC)}$. The case $\ka=0$ yields  homogeneous estimates. 
Eventually, uniqueness again implies that one can extend fundamental solutions to larger  strips.

\section{Stability under perturbation}
\label{sec:perturbations}

 Our methodology allows us to prove a general stability result.
 
\begin{thm}[Stability]\label{thm:stability} Let $\Bdot, \cA$ as in $\mathbf{(B)}$ and $\mathbf{(A)}$ and $\ka\ge c^{a}$, $\ka>0$. Assume that 
\begin{align*}
\| ((\partial_{t}+v\cdot \nabla_{x}) +\cA+ \ka^2)^{-1}\|_{\Wbk\to \Vbk} <\infty.
\end{align*}
Let $\Pdot$ be a  space as in $\mathbf{(B)}$ and $\cP:\Pdot \to \Pdot'$ a continuous operator built as in $\mathbf{(A)}$ with form $p_{t,x}$, and \eqref{eq:bounda} with $c^{p}=0$ and  $\Lambda^{p}=\|\cP\|_{\Pdot \to \Pdot'}$. If $\Lambda^{p}$ is small enough, then 
\begin{align*}
\| ((\partial_{t}+v\cdot \nabla_{x}) +\cA+ \cP+ \ka^2)^{-1}\|_{\Wbk\to \Vbk} \le M\| ((\partial_{t}+v\cdot \nabla_{x}) +\cA+ \ka^2)^{-1}\|_{\Wbk\to \Vbk},
\end{align*}
where $M$ depends only on $d,\beta,\Lambda^a$ and embedding constants. In particular, if the bound is uniform in $\ka$ in the assumption it is also uniform for the conclusion.

If $(\partial_{t}+v\cdot \nabla_{x}) +\cA+ \ka^2$ is causal, so is $(\partial_{t}+v\cdot \nabla_{x}) +\cA+ \cP+ \ka^2$. 

In particular, if  $\mathbf{(WPC)}$ holds for $\cA$ with some $\ka>0$, and $\cP$ has small enough norm, then $(\partial_{t}+v\cdot \nabla_{x}) +\cA+ \cP$ has a unique fundamental solution. 
 \end{thm}
 
\begin{rem}
 The restriction $\ka>0$ is not needed when $\beta<d/2$ or when $\cD(\Omega)\subset \Bdot'$.  See Section~\ref{sec:duality}. 
\end{rem}

\begin{proof} First, we have 
$\cA+\cP+\ka^2 : \E_{\ka}^{\vphantom{'}}=\Bka^{\vphantom{'}}\cap\,\Pdot \to \E_{\ka}'=\Bka' + \,\Pdot'$.   The equation
\[ ((\partial_{t}+v\cdot \nabla_{x}) +\cA+ \cP+ \ka^2)f=S\] can be posed for weak solutions in $\E_{\ka}$.
As $\cK_{\cA+\ka^2}:=((\partial_{t}+v\cdot \nabla_{x}) +\cA+  \ka^2)^{-1}: \Wbk \to \Vbk$, we have  by the embeddings
$\cK_{\cA+\ka^2}: \E_{\ka}' \to \E_{\ka}^{\vphantom{'}}$. Call 
 its bound $C_{0}$.   
Now  the perturbation $\cP$ can be seen as an operator from $\E_{\ka}^{\vphantom{'}}$ to $\E_{\ka}'$ with bound not exceeding $\|\cP\|_{\Pdot\to \Pdot'}$. If its norm is bounded by  $(2C_{0})^{-1}$,   
 $\|\cK_{\cA+\ka^2} \cP\|_{\E_{\ka}\to \E_{\ka}} \le 1/2$ and  $\| \cP\cK_{\cA+\ka^2}\|_{\E_{\ka}'\to \E_{\ka}'} \le 1/2$. We claim that  
 \begin{equation}
\label{eq:firstbddness}
((\partial_{t}+v\cdot \nabla_{x}) +\cA+\cP+ \ka^2)^{-1}= (1+\cK_{\cA+\ka^2} \cP)^{-1}\cK_{\cA+\ka^2}: \E_{\ka}' \to \E_{\ka}^{\vphantom{'}}
\end{equation}
with bound $2C_{0}$. 
The bound follows from the equality once we show that the right hand side yields a unique solution to the corresponding equation.  

For the existence, let $S\in \E_{\ka}'$. We know that $\cK_{\cA+\ka^2} S\in \Vbk\subset \E_{\ka}^{\vphantom{'}}$ by Theorem~\ref{thm:optimalregularity}.  Expand in a Neumann series  $f=(1+\cK_{\cA+\ka^2}  \cP)^{-1}\cK_{\cA+\ka^2} S$ and observe  that 
the $n$th partial sum $f_{n}\in \E_{\ka}$ is a weak solution to  \[((\partial_{t}+v\cdot \nabla_{x}) +\cA+ \cP+\ka^2)f_{n}= S + (- \cP \cK_{\cA+\ka^2} )^{n+1}S.\] Hence, taking limits  $n\to \infty$ in the weak formulation using  $\|(- \cP \cK_{\cA+\ka^2} )^{n+1}S\|_{\E_{\ka}'} \to 0$ and $\cD(\Omega)\subset \E_{\ka}'\subset \cD'(\Omega)$ as $\ka>0$, we find that $f$ is a weak solution to  $((\partial_{t}+v\cdot \nabla_{x}) +\cA+\cP + \ka^2)f= S$.

For uniqueness, we proceed as follows. Let $f\in \E_{\ka}$ be a weak solution of $((\partial_{t}+v\cdot \nabla_{x}) +\cA+\cP + \ka^2)f= 0$. Let $S\in \E_{\ka}'$. We want to show $\langle f, S \rangle=0$. By Proposition~\ref{prop:duality}, we have that $(-(\partial_{t}+v\cdot \nabla_{x}) +\cA^*+  \ka^2)^{-1}: \Wbk \to \Vbk$ and from what precedes, there exists $g\in \E_{\ka}$ such that $-(\partial_{t}+v\cdot \nabla_{x})g +(\cA+\cP + \ka^2)^*g= S$. As in the proof of Proposition~\ref{prop:duality}, we use the absolute continuity of $t\mapsto \langle  f(t), g(t)\rangle _{\L^2_{x,v}}$ to conclude. 

We have shown \eqref{eq:firstbddness} and applying Proposition~\ref{prop:duality} as $\ka>0$, we obtain the desired conclusion, which  preserves uniformity.  

\medskip

Stability of causality can be done in two steps. We first assume $S\in \E_{\ka}'$ and by \eqref{eq:firstbddness} the weak solution  $f$ is built via the converging Neumann series $$f=(1+\cK_{\cA+\ka^2}\cP )^{-1}\cK_{\cA+\ka^2} S= \sum_{\ell=0}^\infty (-\cK_{\cA+\ka^2}\cP )^{\ell}\cK_{\cA+\ka^2}S.
$$
By definition  of  $\cP$,  for any interval $\I$, we have $\cP g=0$ on $\I\times \,\R^{2d}$ if $g=0$ on $\I\times \,\R^{2d}$. 
By the causality assumption, if $h=0$ on $(-\infty,T)\times \, \R^{2d}$  then  
 $\cK_{\cA+\ka^2}h$ also vanishes there,  thus so does $f$ by induction and summation. 
 
 In the second step, we assume that $S\in \Wbk$. One can approximate $S$ by a sequence  of elements  $S_{j}\in \cS(\Omega)$ which vanish on $\I_{j}\times \,\R^{2d}$ where $\I_{j}=(-\infty, T-1/j)$ if $\I=(-\infty, T)$. The associated  weak solutions $f_{j}$ vanish on $\I_{j}\times \,\R^{2d}$ for each $j$ by the previous step. 
We have convergence  $f_{j}\to f$ in $\Vbk$, so it is easy to conclude that $f$ vanishes on $\I\times \,\R^{2d}$.

\

We have shown that  $\mathbf{(WPC)}$ holds for $\cA+\cP$ with the same $\ka$ as for $\cA$. By Theorem~\ref{thm:fso},  $(\partial_{t}+v\cdot \nabla_{x}) +\cA+\cP$ has a  fundamental solution.
\end{proof}

\section{Well-posedness under appropriate assumptions}

We have developed a theory to obtain a fundamental solution for $(\partial_{t}+v\cdot \nabla_{x}) +\cA$ under assumption $\mathbf{(WPC)}$  and it is stable under perturbation. The next step is to construct $\cA$ for which $\mathbf{(WPC)}$ can be checked. 
We introduce two general conditions.

\

 The first condition is the case of operators with bounded coefficients considered in~\cite{AIN}. 
  
  \

\paragraph{\textbf{Hypothesis (H1): Coercive operators with bounded coefficients}}  Consider $\Bdot= \Fbd$ in $\mathbf{(B)}$  and  $\cA$ in $\mathbf{(A)}$ such that  
 for some $0<\lambda_{a}< \Lambda^{a}<\infty$,  $c_{a},c^{a}\in [0,\infty)$, and  uniformly for almost all $(t,x)\in \R\times \R^d$, 
 \begin{equation}
\label{e:ellip-upper-inhom}
\forall f,g \in \Hb_{v}, \quad  |a_{t,x}(f,g)|\le \Lambda^{a} (\|f\|_{\Hbd_{v}}+c^{a}\|f\|_{\L^2_{v}})( \|g\|_{\Hbd_{v}}+c^{a}\|g\|_{\L^2_{v}})
\end{equation}
and 
\begin{align}
\label{eq:weakelliptic-bis}
 \forall f \in \Hb_{v}, \quad  \Re a_{t,x}(f,f) \ge \lambda_{a} \|f\|_{\Hdot^\beta_{v}}^2 -c_{a}^2\|f\|^2_{\L^2_{v}}.
 \end{align}
Note that \eqref{e:ellip-upper-inhom} is stronger than \eqref{eq:bounda}.

\

The next condition is a variant of $\mathbf{(H1)}$ where operators may not have bounded coefficients.   

\

\paragraph{\textbf{Hypothesis (H2): Coercive operators with form lower bound}}
We assume that $\Bd$ in $\mathbf{(B)}$ is a reflexive Banach space whose norm is equivalent to the $\Fbd$ norm  and  let $\cA$ defined in $\mathbf{(A)}$ such that  \eqref{e:ellip-upper-inhom} and 
\eqref{eq:weakelliptic-bis} hold only for $f,g\in \B_{v}$.

\begin{rem}  $\mathbf{(H1)}$ is the case considered in 
  \cite{AIN}. Moreover,  if the constants $c_{a}$ and $c^{a}$ both vanish, this  is a homogeneous situation for $\cA$ also considered in \cite{AIN}. 
\end{rem}
 
 \begin{rem}  In $\mathbf{(H2)}$, the space $\Bdot_{v}$ has a norm equivalent to that of $\Hdot^\beta_{v}$. 
\end{rem}

The following  result  allows  to implement our theory. It includes both well-posed\-ness and the best possible regularity result in this generality.

\begin{thm}[Well-posedness and optimal regularity] \label{thm:exuniq} Assume $\mathbf{(H1)}$ or $\mathbf{(H2)}$. Then $\mathbf{(WPC)}$ holds for all $ \ka\ge \ka_{a}:= \max (2c_{a},c^a)$ and  there exists 
 $0<C<\infty$ depending on $d,\beta, \Lambda^a, \lambda_{a}, c_{1},c_{2}$  such that for all $\ka\ge \ka_{a}$,    $((\partial_{t}+v\cdot \nabla_{x}) +\cA+ \ka^2)^{-1}: \Wbk  \to \Vbk$ with bound $C$. In particular, $(\partial_{t}+v\cdot \nabla_{x}) +\cA$ has a fundamental solution.

The same holds for 
 $(-(\partial_{t}+v\cdot \nabla_{x}) +\cA^*+ \ka^{2})^{-1}$ with causality changed to anti-causality. 
  \end{thm}

\begin{rem}
 It is not clear whether we can take $c_{a}$ instead of $2c_{a}$ in the definition of $\ka_{a}$. Taking a multiple of $c_{a}$ makes $\ka_{a}=0$ a possible value. 
\end{rem}
 
\begin{proof} We proceed under $\mathbf{(H1)}$, then $\mathbf{(H2})$, to show $\mathbf{(WPC)}$ and obtain the uniform bound in $\ka$. 

 \medskip
 
\paragraph{\textbf{Case  1: Assume $\mathbf{(H1)}$}}

\subparagraph{\emph{Step 1 (uniqueness in $\Fbk$ and causality):}}   Assume $f$ is a weak solution in $\Fbk$ with $S=0$. Fix $-\infty<\sigma\le \tau<+\infty$. The energy equality between $\sigma\le \tau$ yields
 \[
   \|f(\tau)\|_{\L^2_{x,v}}^2 + 2 \int_{\sigma}^\tau \int_{\R^d} \Re(a_{t,x}(f,f)+\ka^2\angle {f}{f})\, \dx\dt =  \|f(\sigma)\|_{\L^2_{x,v}}^2,
 \]
 writing $f$ for $f(t,x,\cdot)$ in the integrand. As $\ka \ge c_{a}$, the integrand is nonnegative by \eqref{eq:weakelliptic-bis}, hence $\|f(\tau)\|_{\L^2_{x,v}}^2 \le  \|f(\sigma)\|_{\L^2_{x,v}}^2$. 
Using the 0 limit at $\sigma\to - \infty$ in $\L^2_{x,v}$ shows that $\|f(\tau)\|_{\L^2_{x,v}}^2=0$. Hence $f=0$ identically.

 For the causality, we assume $S=0$ on $(-\infty, T)\times \, \R^{2d}$ and $S\in \Wbk$. The energy equality yields the same identity for $\tau<T$.  
 Hence, $f=0$ on $(-\infty, T)\times \, \R^{2d}$. 

\medskip

\subparagraph{\emph{Step 2 (existence for $(\Fbk)'$ sources):}}   If $\ka=0$ (hence $\ka_{a}=0$) this is already in \cite{AIN}, Corollary~1.11.    
  We assume from now on  $\ka \ge 2c_{a}$ and $\ka>0$. Again, the result is essentially in \cite{AIN}.  However, we give the argument to make  the  bounds  explicit and  uniform in $\ka$ so that one can pass to the limit when $\ka\to 0$ on applying Proposition~\ref{prop:limit} (This is the argument in \cite{AIN} when $\beta\ge d/2$.)
  
  \smallskip
We apply the Lions representation theorem on the  space $\Fbk$ (see the recent short proof in \cite{ACE} based on the Hahn-Banach theorem), renormed with $\|f\|_{\Fbk}=(\|f\|_{\Fbd}^2+ \ka^2 \|f\|_{\L^2_{\vphantom{t,x} t,x,v}}^2)^{1/2}$ to make it a Hilbert space. 
 The sesquilinear form on $\Fbk \times \cD(\Omega)$,  
 $$
 E(f,h)= -\langle  f , (\partial_{t}+v\cdot\nabla_{x})  h\rangle + a(f,h) +\ka^2 \angle{f}{h},
 $$
 is continuous for fixed $h$ as a function of $f$ and satisfies for $h\in \cD(\Omega)$
 \begin{align*}
 \Re E(h,h) \ge \lambda_{a} \|h\|_{\Fbd}^2 +(\ka^2-c_{a}^2)\|h\|_{\L^2_{\vphantom{t,x} t,x,v}}^2 
 \ge \min(\lambda_{a}, 3/4) \|h\|_{\Fbk}^2
 \end{align*}
 since $\ka^2-c_{a}^2\ge3\ka^2/4$.
 Thus for  $S\in (\Fbk)'$   
 there exists $f\in \Fbk$ such that $E(f,h)=\angle{S}{h}$   for all $h\in \cD(\Omega)$ and  $$\|f\|_{\Fbk} \le \min(\lambda_{a}, 3/4)^{-1} 
 \|S\|_{(\Fbk)'}.
 $$

 \smallskip
 
\subparagraph{\emph{Step 3 (conclusion):}} Step 2 proves that  $$((\partial_{t}+v\cdot \nabla_{x}) +\cA+ \ka^2)^{-1}:(\Fbk)' \to \Fbk$$ with uniform estimate if $\ka \ge 2c_{a}$ and $\ka>0$.   As we assume   $\ka \ge c^a$ as well,  we can apply Proposition~\ref{prop:duality}.  We get the conclusion, also with uniform estimate if $\ka \ge \ka_{a}$ and $\ka>0$.
 
 \medskip

\paragraph{\textbf{Case 2: Assume $\mathbf{(H2)}$}} Uniqueness in $\Bka$  and causality are proved as above using that the ${\Bka}$-norm  is equivalent to the  $\Fbk$-norm. For existence, we apply again
the Lions representation  \cite[Theorem~2.1]{ACE} with $\X=\Bka$ which is a reflexive Banach space and  $\Y=\cD(\Omega)$ equipped with the ${\Bka}$-norm. The sesquilinear form on $\Bka \times\,  \cD(\Omega)$,  
 $$
 E(f,h)= -\langle  f , (\partial_{t}+v\cdot\nabla_{x})  h\rangle + a(f,h) +\ka^2 \angle{f}{h},
 $$
 is continuous from $\Bka$ to $\IC$ for fixed $h\in \cD(\Omega)$ as a function of $f$ because $f\in \Bka\subset \cD'(\Omega)$, $g=(\partial_{t}+v\cdot\nabla_{x})  h\in \cD(\Omega)$ and $(\cA+\ka^2)f \in \Bka' \subset \cD'(\Omega)$. Next, for
 $h\in \cD(\Omega)$
 \begin{align*}
 \Re E(h,h) \ge \lambda_{a} \|h\|_{\Fbd}^2 +(\ka^2-c_{a}^2)\|h\|_{\L^2_{\vphantom{t,x} t,x,v}}^2 
 \ge \min(\lambda_{a}, 3/4) \|h\|_{\Fbk}^2 \sim \min(\lambda_{a}, 3/4) \|h\|_{\Bka^{\vphantom{\beta}}}^2
 \end{align*}
 when $\ka\ge 2 c_{a}$. Since $\cD(\Omega)\subset \B_{\ka}$, this implies $\sup\{|E(f,h)| \, ;\,  \|f\|_{\Bka}\le1\} \ge  c \|h\|_{\Bka}$ for some $c>0$.  
 Hence,  for  $S\in \Bka'$   
 there exists $f\in \Bka$ such that $E(f,h)=\angle{S}{h}$ with $\|f\|_{\Bka}^{\vphantom{'}}\le c^{-1}\|S\|_{\Bka'}$.
 Thus, we have proved 
 $$((\partial_{t}+v\cdot \nabla_{x}) +\cA+ \ka^2)^{-1}:\Bka'\to \Bka^{\vphantom{'}}$$ 
 with uniform bound  when $\ka\ge 2c_{a}$. 
 From here, the conclusion is the same as in Case~1.
 \end{proof}

  \section{Examples}
  \label{sec:examples}

  We now come to examples, illustrating the advantage of having adopted such an abstract framework. At this level of generality and given the required embeddings in $\mathbf{(B)}$, a natural idea is to assume lower order terms in  mixed Lebesgue spaces and even mixed Lebesgue-Lorentz spaces for lower order coefficients of local diffusions and lower order kernels of non-local diffusions. This is based on anisotropic Sobolev embeddings thanks to O'Neil's results for Riesz potentials \cite{O'Neil1963convolution}, see also \cite{Tartar98} and \cite{PP} for related results. 
  
  Although we do not give details in absence of a specific application, it is possible to consider higher order $\beta>  1$ subdiffusive operators and apply the strategy developed here as we only assumed $0<\beta$ in our theory.  

\subsection{Local diffusions}  

The local diffusions when $\beta=1$  are given by 
  \[ 
\cA= -\div_{v}A \nabla_{v} -\div_{v}\oa  +\ob\cdot  \nabla_{v}+\oc. 
\]
We assume that $A=A(t,x,v)$ is a $d\times d$-matrix  with bounded, possibly complex, measurable coefficients. We set $\|A\|_{\infty}=\Lambda_{A}$ and we assume coercivity with lower bound $\lambda_{A}>0$, that is $\Re\langle A\xi,\xi \rangle \ge \lambda_{A}|\xi|^2$ for all $\xi\in \IC^d$.  The vector fields  $\oa$, $\ob$ and the potential $\oc$ have measurable, possibly  complex, entries. We  use when needed
that $\|f\|_{\L^2_{t,x}\Hdot^1_{\vphantom{t}v}}\sim \|\nabla_{v}f\|_{\L^2_{t,x,v}}$.  Unless specified,  $d\ge 1$.

  \subsubsection{Example 1}

We assume   
\[
|\oa|^2, |\ob |^2, \oc \in (\L^{\tilde r}_{t}\L^{\tilde p,\infty}_{x}\L^{\tilde q, \infty}_{v})_{\varepsilon} + \L^\infty_{t,x,v} \ 
\]
with  the scaling-critical relation
\[
\frac{2}{\tilde r} +\frac{3d}{\tilde p} +\frac{d}{\tilde q}= 2, \quad  \tilde r, \tilde p,\tilde q\in (1,\infty].
\]
and small enough $\varepsilon>0$.   More precisely, we assume the following condition, explaining the $\varepsilon$ in our notation.\\

\paragraph{\textbf{Assumption $\mathbf{(D)}$:}}  The coefficients $\oa, \ob, \oc$ can be decomposed, for  $\varepsilon>0$ small enough, as  
\[
\oa=\oa_{\varepsilon}+\oa_{\infty}, \ \ob=\ob_{\varepsilon}+\ob_{\infty}, \ \oc=\oc_{\varepsilon}+\oc_{\infty}
\]
with 
\[
P_{\varepsilon}:=\||\oa_{\varepsilon}|^2\|^{1/2}_{\L^{\vphantom{\tilde p}\tilde r}_{t}\L^{\tilde p,\infty}_{\vphantom{t} x}\L^{\tilde q, \infty}_{\vphantom{t} v}}+ \||\ob_{\varepsilon}|^2\|^{1/2}_{\L^{\vphantom{\tilde p}\tilde r}_{t}\L^{\tilde p,\infty}_{\vphantom{t}x}\L^{\tilde q, \infty}_{\vphantom{t}v}}+ \|\oc_{\varepsilon}\|_{\L^{\vphantom{\tilde p}\tilde r}_{t}\L^{\tilde p,\infty}_{\vphantom{t}x}\L^{\tilde q, \infty}_{\vphantom{t}v}} \le \varepsilon
\]
and 
\[
 \|\oa_{\infty}\|_{\infty}+\|\ob_{\infty}\|_{\infty}+ \|\oc_{\infty}\|_{\infty}<\infty. 
\]\\

The space $\L^{p,q}$ is the Lorentz space on $\R^d$
 defined as the set of all measurable functions $f$ on $\mathbb{R}^d$ such that $\left\| f\right\|_{\L^{p,q}} < \infty$, with
\begin{align*}
\left\| f\right\|_{\L^{p,q}} =\begin{cases}
\left (q \int_{0}^{\infty} s^{q-1} \lambda_f(s)^{q/p} \, \mathrm{d}s \right )^{1/q}, &   p,q \in [1,\infty)  , \\ 
 \sup_{s>0} s \, \lambda_f(s)^{1/p}, & p \in [1,\infty), \  q=\infty,
 \\ 
 \|f\|_{\L^\infty}, &  p=\infty, q \in [1,\infty], 
 \end{cases} 
 \end{align*}
 where 
$$ \lambda_f(s)=\mu\left ( \left\{ x\in \mathbb{R}^d : |f(x)|>s\right\} \right ), 
$$ and $\mu$ is here Lebesgue measure. It is known that $\L^{p,p}=\L^{p}$ isometrically and $\left\| f\right\|_{\L^{p,q}}$ is non-increasing as a function of $q$, so that $\L^{p,q_1}\subset \L^{p,q_2}$ if $1\le q_1 
\leq q_2\le \infty$. When $1<p,q<\infty$, $(\L^{p,q})'=\L^{p',q'}$.   We refer to \cite[Ch. V, \S3 ]{SW} for  details.

\begin{prop}[Typical example]
 If 
\[
|\oa|^2, |\ob |^2, \oc \in \L^{q_{1}}_{t}\L^{q_{2},s_{2}}_{\vphantom{t}x}\L^{q_{3}, s_{3}}_{\vphantom{t}v}
\]
with $q_{1}, q_{2}, q_{3}\in (1,\infty], \ s_{2},s_{3} \in [1,\infty]$ and
either the subcritical relation
\[
\frac{2}{q_{1}} +\frac{3d}{q_{2}} +\frac{d}{q_{3}}< 2, \] 
or the critical relation 
\[
\frac{2}{q_{1}} +\frac{3d}{q_{2}} +\frac{d}{ q_{3}}= 2, \quad \mathrm{all\  five\ exponents    \  finite},
\] 
then $\mathbf{(D)}$ holds.
\end{prop}

\begin{proof}
Truncate the coefficients above some height to produce smallness of $P_{\varepsilon}$.   \end{proof}

If  the critical relation holds and  one of  the exponents is infinite,  then the truncations may not provide smallness.  Here are two particular situations where $\mathbf{(D)}$ holds.   
\begin{enumerate}
\item (\textbf{Critical relation}) Assume  $\tilde r=\tilde p=\tilde q=2d+1$ and consider  
  $|\oa|^2, |\ob |^2, \oc \in \L^{2d+1}_{t,x,v}$, with arbitrary size. This corresponds to the critical relation with $s_{2}=s_{3}=2d+1$ finite and we may use the proposition. 
   
\item (\textbf{Subcritical relation}) Consider $|\oa|^2, |\ob |^2, |\oc| \in \L^{\infty}_{t,x}\L^{\tilde q,\infty}_{v}$ where $\tilde q>d/2$ with $d\ge 2$. This corresponds to a subcrititical relation with  $q_{1}=\infty, q_{2}=\infty$, $q_{3}=\tilde q$. We may use the proposition  to decompose into a term  small in $\L^{\infty}_{t,x}\L^{d/2,\infty}_{v}$  plus  a bounded term.  For example, we can take  each 
$|\oa|^2, |\ob |^2, |\oc|$  bounded by any multiple of the potential $ |v|^{-d/\tilde q}$. 
\end{enumerate}

We next check that Assumption $\mathbf{(D)}$ implies a decomposition 
$\cA=\widetilde \cA+\cP$ where $\widetilde \cA$ satisfies  the hypothesis $\mathbf{(H1)}$ and $\cP$ is a perturbation   with a suitable space $\Pdot$.  

Let $r=2(\tilde r)'$, $p=2(\tilde p)'$, $q=2(\tilde q)'$, where the exponents $\tilde r, \tilde p, \tilde q$ are the ones in $\mathbf{(D)}$ and $'$ denotes the H\"older conjugation. Note that
\[
\frac{2}{r} +\frac{3d}{p} +\frac{d}{q}= 2d, \quad r,p,q\in [2,\infty).
\]
We set 
\[
\Pdot= \L^{\vphantom{p}r}_{t}\L^{p,2}_{ \vphantom{t}x}\L^{q, 2}_{\vphantom{t} v} \cap \L^2_{t,x}\Hdot_{ v}^1.
\]

\begin{prop}\label{prop:H1}
 The space $\Pdot$ satisfies $\mathbf{(B)}$. 
\end{prop}

We need the following lemma.

\begin{lem}\label{lem:ex1} {}\ 
\begin{enumerate}
\item For $p_{0}, q_{0}\in [2,\infty)$ with $\frac{3d}{p_{0}} +\frac{d}{ q_{0}}=2d-1$, 
we have the anisotropic  Sobolev embedding
\[
\L^2_{x}\Hdot_{ v}^1
\cap  \L^2_{v}\Hdot_{\vphantom{x} x}^{{1}/{3}} \subset \L^{p_{0},2}_{\vphantom{t}x}\L^{q_{0},2}_{\vphantom{t}v}. \]
\item For $p,q,r\in [2,\infty)$ with $\frac{2}{r} +\frac{3d}{p} +\frac{d}{q}= 2d$,  we have
 the continuous inclusion
 \[(\Zdot^1)' \subset  \L^{\vphantom{p}r}_{t}\L^{p,2}_{ \vphantom{t}x}\L^{q, 2}_{\vphantom{t} v}.
 \]

\end{enumerate}

\end{lem}

\begin{proof} For (i), we first use \cite{AN}, Lemma 2.6 and Proposition 2.8. If $f$ has Fourier transform $\hat f(\varphi,\xi)$ compactly supported away from $\varphi=0$ and $\xi=0$,   then $f\in \L^2_{x}\Hdot_{ v}^1
\cap  \L^2_{v}\Hdot_{\vphantom{x} x}^{{1}/{3}}$ if and only if the Fourier transform satisfies $\hat f(\varphi,\xi)= (|\varphi|^{1/3}+|\xi|)^{-1} \hat g(\varphi,\xi)$ with $g\in \L^2_{x,v}$. A computation based on the corresponding anisotropic Littlewood-Paley decomposition yields that 
\[
|f(x,v)|\le c_{d}\int_{\R^{2d}} \frac{ |g(y,w)|}{(|x-y|^{1/3}+|v-w|)^{4d-1}}\dd(y,w),
\]
with $c_{d}>0$ a dimensional constant.  Pick $\alpha_{0},\beta_{0}\in (0,1)$ with $\alpha_{0}+\beta_{0}=1$ and
\[
\frac 1 2 -\frac 1 {p_{0}}= \frac {\alpha_{0}}{3d}, \quad \frac 1 2 -\frac 1 {q_{0}}= \frac {\beta_{0}}{d}.
\]
Then we can write 
\[
|f(x,v)|\le c_{d}\int_{\R^{2d}} \frac{ |g(y,w)|}{|x-y|^{d-\alpha_{0}/3}|v-w|^{d-\beta_{0}}}\dd(y,w),
\]
 and thus
 \[|f|\lesssim I_{\alpha_{0}/3, x}I_{\beta_{0},v}|g|,
 \]
 where $I_{\gamma,z}$ is the Riesz potential of order $\gamma\in (0,d)$ in the $z$-variable on $\R^d$. By O'Neil's result $I_{\gamma,z}: \L^2\to \L^{p,2}$ when $\frac 1 2 -\frac 1 {p}= \frac {\gamma}{d}$ \cite{O'Neil1963convolution}. Together with the relation between $p_{0},\alpha_{0}$ and $q_{0},\beta_{0}$, we obtain
 \begin{align}
\label{eq:switch} \| I_{\alpha_{0}/3, x}I_{\beta_{0},v}|g|\|_{\L^{p_{0},2}_{\vphantom{t}x}\L^{q_{0},2}_{\vphantom{t}v}}
 &= \| \|I_{\beta_{0},v} I_{\alpha_{0}/3, x}|g|\|_{\L^{q_{0},2}_{\vphantom{t}v}}\|_{\L^{p_{0},2}_{\vphantom{t}x}}\\
\nonumber & \lesssim  \| \| I_{\alpha_{0}/3, x}|g|\|_{\L^{\vphantom{p_{0}}2}_{\vphantom{t}v}}\|_{\L^{p_{0},2}_{\vphantom{t}x}}\\
 \nonumber &\le  \| \| I_{\alpha_{0}/3, x}|g|\|_{\L^{p_{0},2}_{\vphantom{t}x}} \|_{\vphantom{\L^{p_{0},2}_{\vphantom{t}x}}\L^{\vphantom{p_{0}}2}_{\vphantom{t} v}}
 \lesssim \| \| g\|_{\L^{\vphantom{p_{0}}2}_{\vphantom{t}x}} \|_{\L^{\vphantom{p_{0}}2}_{\vphantom{t}v}}
\end{align}
where we used $p_{0}\ge 2$ in the next to last inequality (Minkowski). We conclude by density. 

To prove (ii), we recall that 
$ (\Zdot^1)'= \L^2_{t,x}\Hdot_{ v}^1
\cap  \L^2_{t,v}\Hdot_{\vphantom{t^a} x}^{{1}/{3}} \cap \L^\infty_{t}\L^2_{\vphantom{t} x,v}$. If $r=2$, this follows directly from (i) with $(p_{0},q_{0})=(p,q)$. Assume next $r>2$. 
Taking $\theta= 2/r\in (0,1)$, we see that
\[
\bigg(\frac 1 r, \frac 1 p,\frac 1 q\bigg)= (1-\theta) \bigg(\frac 1 \infty, \frac 1 2,\frac 1 2\bigg) + \theta \bigg(\frac 1 2, \frac 1 {p_{0}},\frac 1 {q_{0}}\bigg)
\]
for a pair $(p_{0},q_{0})$ as in (i). For almost all $t$, we have 

\[\|f\|_{\L^{p,2}_{\vphantom{t}x}\L^{q,2}_{\vphantom{t}v}} \le \|f\|_{\L^{p_{0},2}_{\vphantom{t}x}\L^{q_{0},2}_{\vphantom{t}v}}^\theta \|f\|_{\L^{\vphantom{p_{0}}2}_{\vphantom{t}x}\L^{\vphantom{p_{0}}2}_{\vphantom{t}v}}^{1-\theta}
\le \|f\|_{\L^{p_{0},2}_{\vphantom{t}x}\L^{q_{0},2}_{\vphantom{t}v}}^\theta \| f\|_{\L^{\vphantom{p_{0}}\infty}_{t}\L^{\vphantom{p_{0}}2}_{\vphantom{t}x,v}}^{1-\theta}.
\]
Using $r\theta=2$, we obtain
\[\|f\|_{\L^{\vphantom{p_{0}}r}_{t}\L^{\vphantom{p_{0}}p,2}_{\vphantom{t} x}\L^{\vphantom{p_{0}}q,2}_{\vphantom{t} v}} \le \|f\|_{\L^{\vphantom{p_{0}}2}_{t}\L^{p_{0},2}_{\vphantom{t} x}\L^{q_{0},2}_{\vphantom{t} v}}^\theta  \| f\|_{\L^{\vphantom{p_{0}}\infty}_{t}\L^{\vphantom{p_{0}}2}_{\vphantom{t}x,v}}^{1-\theta}.
\]
Using (i), the conclusion follows. 
 \end{proof}

\begin{proof}[Proof of Proposition~\ref{prop:H1}]The assertions (i) to (iii) are straightforward. We only consider the embeddings.
The embedding \eqref{eq:B} follows from definition and (ii) in 
Lemma~\ref{lem:ex1}.  Concerning  \eqref{eq:B'}, we start from 
\[
\Pdot'= \L^{\vphantom{p'}r'}_{t}\L^{p',2}_{ \vphantom{t}x}\L^{q', 2}_{\vphantom{t} v} +\L^2_{t,x}\Hdot_{ v}^{-1}.
\]
The  inclusion $\Pdot' \supset \L^2_{t,x}\Hdot_{ v}^{-1}$ is clear. The  inclusion $\Pdot'\subset \Zdot^1$ follows from a  functional analysis argument. If $S\in \Pdot'$ and  $\theta\in \cD(\R^{2d})$, then 
$(S,\theta)$ can be identified to a locally integrable function of $t$. We also have $\Ydot^1\subset (\Zdot^1)' \subset \Pdot$, hence 
$\Bd'\subset  \L^2_{t,x}\Hdot_{ v}^{-1}+ \L^2_{t,v}\Hdot_{ x}^{-1/3}+ \cM(\L^2_{x,v})$. Writing $S=S_{1}+S_{2}+S_{3}$ according to this sum, 
we see that the term $(S_{3},\theta)$ is a bounded measure on $\R$ and also a locally integrable function, so it is in fact an integrable function on $\R$. The desired inclusion follows easily from this. 
 \end{proof}

\begin{prop}\label{lem:H2}
 Under the condition $\mathbf{(D)}$, the operator $\cA$ is a perturbation of an operator satisfying $\mathbf{(H1)}$.
\end{prop}
 
\begin{proof}
 Fix $\varepsilon>0$. We decompose the coefficients according to $\mathbf{(D)}$. We use that the Lorentz spaces obey the same rules as Lebesgue  spaces for  H\"older's inequality: $$\L^{p_{1},q_{1}}\times 
\L^{p_{2},q_{2}}\times \L^{p_{3},q_{3}}\subset \L^1, \quad \frac 1 {p_{1}}+
\frac 1 {p_{2}}+\frac 1 {p_{3}}=1,\quad \frac 1 {q_{1}}+
\frac 1 {q_{2}}+\frac 1 {q_{3}}=1.$$ We also observe that $\|\nabla_{v}f\|_{\L^2_{t,x,v}}= \|f\|_{\L^2_{t,x}\Hdot_{\vphantom{t} v}^1}$. Hence, with the definition of the exponents, 
we see that for $f, g\in \Pdot$, 
\begin{align*}
\int_{\R\times \,\R^{2d}} |\oa_{\varepsilon} f \cdot \nabla_{v}g | \dd(t,x,v)&\le 
\|\oa_{\varepsilon}\|_{\L^{\vphantom{\tilde p}2 \tilde r}_{t}\L^{2\tilde p,\infty}_{\vphantom{t} x}\L^{2\tilde q, \infty}_{\vphantom{t} v}}
\|f\|_{\L^{\vphantom{p}r}_{t}\L^{p,2}_{ \vphantom{t}x}\L^{q, 2}_{\vphantom{t} v}} \|g\|_{\L^2_{t,x}\Hdot_{\vphantom{t}  v}^1} \\
&\le\||\oa_{\varepsilon}|^2\|^{1/2}_{\L^{\vphantom{\tilde p}\tilde r}_{t}\L^{\tilde p,\infty}_{\vphantom{t} x}\L^{\tilde q, \infty}_{\vphantom{t} v}} 
\|f\|^{\vphantom{1/2}}_{\Pdot^{\vphantom{\tilde p}}_{\vphantom{t}}}
\|g\|^{\vphantom{1/2}}_{\Pdot^{\vphantom{\tilde p}}_{\vphantom{t}}}\\
\int_{\R\times \,\R^{2d}} |\ob_{\varepsilon} g \cdot \nabla_{v}f | \dd(t,x,v)&\le 
\|\ob_{\varepsilon}\|_{\L^{\vphantom{\tilde p}2\tilde r}_{t}\L^{2\tilde p,\infty}_{\vphantom{t} x}\L^{2\tilde q, \infty}_{\vphantom{t} v}}\|g\|_{\L^{\vphantom{p} r}_{t}\L^{p,2}_{ \vphantom{t}x}\L^{q, 2}_{\vphantom{t} v}} \|f\|_{\L^{2}_{ t,x}\Hdot_{ v}^1} \\
&\le \||\ob_{\varepsilon}|^2\|^{1/2}_{\L^{\vphantom{\tilde p}\tilde r}_{t}\L^{\tilde p,\infty}_{\vphantom{t} x}\L^{\tilde q, \infty}_{\vphantom{t} v}} 
\|f\|^{\vphantom{1/2}}_{\Pdot^{\vphantom{\tilde p}}_{\vphantom{t}}}
\|g\|^{\vphantom{1/2}}_{\Pdot^{\vphantom{\tilde p}}_{\vphantom{t}}}\\
\int_{\R\times \,\R^{2d}} |\oc_{\varepsilon} f g | \dd(t,x,v)&\le 
\|\oc_{\varepsilon}\|_{\L^{\vphantom{\tilde p}\tilde r}_{t}\L^{\tilde p,\infty}_{\vphantom{t} x}\L^{\tilde q, \infty}_{\vphantom{t} v}}\|f\|_{\L^{\vphantom{p} r}_{t}\L^{p,2}_{ \vphantom{t}x}\L^{q, 2}_{\vphantom{t} v}} \|g\|_{\L^{\vphantom{p}r}_{t}\L^{p,2}_{ \vphantom{t}x}\L^{q, 2}_{\vphantom{t} v}} \\
&\le \|\oc_{\varepsilon}\|_{\L^{\vphantom{\tilde p}\tilde r}_{t}\L^{\tilde p,\infty}_{\vphantom{t} x}\L^{\tilde q, \infty}_{\vphantom{t} v}} \|f\|_{\Pdot^{\vphantom{\tilde p}}_{\vphantom{t}}}\|g\|_{\Pdot^{\vphantom{\tilde p}}_{\vphantom{t}}}.
 \end{align*}
 Hence with  $\cP:=-\div_{v}\oa_{\varepsilon}  +\ob_{\varepsilon}\cdot  \nabla_{v}+\oc_{\varepsilon}$ and $P_{\varepsilon}\le \varepsilon$ in $\mathbf{(D)}$,  we obtain  \eqref{eq:bounda} with $c^{p}=0$ and $\Lambda^{p}=\varepsilon$. 
 For the other term, 
 $\widetilde \cA:= -\div_{v}A \nabla_{v} -\div_{v}\oa_{\infty}  +\ob_{\infty}\cdot  \nabla_{v}+\oc_{\infty}$, it is easy to show \eqref{e:ellip-upper-inhom} and \eqref{eq:weakelliptic-bis} with constants $\Lambda^{a}=\Lambda_{A}$, $\lambda_{a}=\lambda_{A}/2$ and appropriate $c^{a},c_{a}$, depending on the $\L^\infty_{t,x,v}$ norms of $\oa_{\infty}, \ob_{\infty}, \oc_{\infty}$. We can thus write $\cA=\widetilde \cA+\cP$ as announced.  
 \end{proof}
 
 According to Remark~\ref{rem:homogeneous}, we can reach homogeneous estimates for critical perturbations with small norm  when working on infinite intervals.  For example, with $\I=\R$,  if 
 $A$ is as above and \[
|\oa|^2, |\ob |^2, \oc \in \L^{\tilde r}_{t}\L^{\tilde p,\infty}_{x}\L^{\tilde q, \infty}_{v}    
\]
with  small enough norms under the scaling-critical relation
\[
\frac{2}{\tilde r} +\frac{3d}{\tilde p} +\frac{d}{\tilde q}= 2, \quad  \tilde r, \tilde p,\tilde q\in (1,\infty],
\]
then 
\[(\pm(\partial_{t}+v\cdot \nabla_{x}) -\div_{v}A \nabla_{v} -\div_{v}\oa  +\ob\cdot  \nabla_{v}+\oc)^{-1} : \Zdot^1\to \Ydot^1
\]
with 
\[  \Zdot^1= \L^2_{t,x}\Hdot^{-1}_{\vphantom{t,x} v} + \L^2_{t,v}\Hdot_{\vphantom{t} x}^{-1/3} +\L^1_{t} \L^2_{x,v}\quad , \quad  \Ydot^1= \L^2_{t,x}\Hdot^{1}_{\vphantom{t,x} v} \cap  \L^2_{t,v}\Hdot_{\vphantom{t} x}^{1/3} \cap \C^{}_{0}\L^2_{x,v}.
\]
For example, take $\oa, \ob$ small multiples of  $|v|^{-2}v$ and $\oc$  small multiple of the Coulomb potential $ |v|^{-2}$. This corresponds to     $\tilde r=\infty, \tilde p=\infty$,  $\tilde q=d/2$ and $d\ge 3$.  
 
\begin{rem}
 We have chosen the same exponents for all coefficients $\oa, \ob,\oc$ but they can be different with the same relations. In that case, the space $\Pdot$ has to be adjusted by taking the intersection of the different $ \L^{\vphantom{p}r}_{t}\L^{p,2}_{ \vphantom{t}x}\L^{q, 2}_{\vphantom{t} v}$ spaces.  
\end{rem}

\subsubsection{Example 2}

 Assume again  $\tilde r=\infty, \tilde p=\infty$,  $\tilde q=d/2$, $d\ge 3$ and $|\oa|^2, |\ob |^2, |\oc| \in \L^{\infty}_{t,x}\L^{d/2,\infty}_{v} + \L^\infty_{t,x,v}$. The decomposition $\mathbf{(D)}$ does not apply because this corresponds to the critical relation with some  infinite exponents. But smallness is not necessary provided we have form lower bounds.  Let $\Bdot= \L^2_{t,x}\L^{\frac{2d}{d-2},2}_{\vphantom{t}v} \cap \L^2_{t,x}\Hdot^1_{\vphantom{t}v}$. Notice that by Sobolev embedding proved by O'Neil in the $v$-variable   $\|f\|_{\Bdot}\sim \|f\|_{\L^2_{t,x}\Hdot^1_{\vphantom{t}v}}$ for $f\in \Bdot$. Moreover,  $\Bdot$ is a reflexive Banach space. Thus we can apply $\mathbf{(H2)}$
 provided that for all $f\in \B^{\vphantom{2}}_{v}=\L^{\frac{2d}{d-2},2}_{\vphantom{t}v} \cap\, \Hdot^1_{\vphantom{t}v}$,  
\[\Re   p_{t,x} (f , f )  \ge -C_{0} \|f\|_{\Hdot^1_{v}}^2 - c_{0}^2\|f\|_{\L^2_{v}}^2 
\]
with $C_{0}<\lambda_{A}$ uniformly in $(t,x)$, where $p_{t,x}$ denotes  the form associated in the $v$-variable to  $-\div_{v}\oa  +\ob\cdot  \nabla_{v}+\oc$. Note that for $f, g\in \B$  
we have 
\[|p_{t,x} (f , g )| \le C^0 (\|f\|_{\L^2_{t,x}\Hdot^1_{\vphantom{t}v}}+ c^0\|f\|_{\L^2_{t,x,v}}) (\|g\|_{\L^2_{t,x}\Hdot^1_{\vphantom{t}v}}+ c^0\|g\|_{\L^2_{t,x,v}}).
\]
For instance, the lower bound can be obtained if $\oa, \ob $ satisfy the decomposition as in  $\mathbf{(D)}$ and $\Re \oc \ge 0$.  Let us isolate another example where this positivity is not necessary.

\begin{cor}[kinetic Schr\"odinger operator]\label{cor:schrodinger} Assume $d\ge 3$. The operator 
\[(\partial_{t}+v\cdot\nabla_{x}) -\Delta_{v}+ c|v|^{-2}\]
 with $c$ a constant or even a bounded function of $t,x,v$ satisfying
$\essinf \Re c> -\frac{(d-2)^2}4$ has a unique fundamental solution $\Gamma(t,s)$ defined and contractive on $\L^2_{x,v}$ for $-\infty<s\le t <\infty$.  
 \end{cor}

\begin{proof}
 Setting $\cA=-\Delta_{v}+ c|v|^{-2}$, we have $\mathbf{(H2)}$ with $\Bdot= \L^2_{t,x}\L^{\frac{2d}{d-2},2}_{\vphantom{t^a} v} \cap \L^2_{t,x}\Hdot^1_{\vphantom{t}v}$ because of Hardy's inequality in the $v$-variable. The contractivity follows from the absolute continuity of $t\mapsto \|\Gamma(t,s)\psi\|^2_{\L^2_{x,v}}$, $t>s$, when $\psi\in \L^2_{x,v}$  and the computation with $u(t)=\Gamma(t,s)\psi$, 
 \[
 \frac {d}{dt} \|u(t)\|^2_{\L^2_{x,v}}= -2 \Re \int_{\R^{2d}} |\nabla_{v}u(t,x,v)|^2+ c(t,x,v)|v|^{-2}|u(t,x,v)|^2\dd(x,v)  \le 0, \quad \text{a.e.}
 \]
 using Fubini's theorem and Hardy's inequality in the $v$-variable.  
\end{proof}

 \subsubsection{Example 3}

This is apparently the same example as in Example 1 but   
with  the condition 
\[
|\oa|^2, |\ob |^2, \oc \in (\L^{\tilde r}_{t}\L^{\tilde q, \infty}_{v}\L^{\tilde p,\infty}_{x})_{\varepsilon} + \L^\infty_{t,x,v}  
\]
under the  critical relation 
\[
\frac{2}{\tilde r} +\frac{3d}{\tilde p} +\frac{d}{\tilde q}= 2, \quad  \tilde r, \tilde p,\tilde q\in (1,\infty]
\]
and $\varepsilon>0$ small enough. 

Here, this is similar to $\mathbf{(D)}$ but we have switched the order of the  Lorentz spaces in the variables $x$ and $v$  while keeping the critical relation. This allows some more critical potentials and drifts. For example, with the exponents  $\tilde r=\infty, \tilde q=\infty$, $\tilde p=3d/2$,  we can take  each 
$|\oa|^2, |\ob |^2, |\oc|$  bounded by $ c|x|^{-2/3}$ with small  constant $c$. 
It is the proof of Lemma~\ref{lem:ex1}, (i), that allows this switch. More precisely, we can change the order of the Lorentz norms in the proof of \eqref{eq:switch} and obtain
\[
\L^2_{x}\Hdot_{ v}^1
\cap  \L^2_{v}\Hdot_{\vphantom{x} x}^{{1}/{3}} \subset  \L^{q_{0},2}_{\vphantom{t}v}\L^{p_{0},2}_{\vphantom{t}x}.
\]
 The space $\Pdot$ is chosen as 
\[
\Pdot= \L^{\vphantom{p}r}_{t}\L^{q,2}_{ \vphantom{t}v}\L^{p, 2}_{\vphantom{t} x} \cap \L^2_{t,x}\Hdot_{\vphantom{t} v}^1.
\]
 We skip further details. 

\begin{rem}
 Again we can choose  different $\tilde p,\tilde q,\tilde r$ for $|\oa|^2, |\ob |^2, |\oc|$. We can also mix with conditions in Examples 1 and 2. The principle is to take the  sum of all spaces  for the coefficients and the intersection of all corresponding solution  spaces.  
\end{rem}

\begin{cor} Assume $d\ge 3$. The perturbed kinetic Schr\"odinger  operator 
\[(\partial_{t}+v\cdot\nabla_{x}) -\Delta_{v}+  c|v|^{-2}+m |x|^{-2/3} \]
  with $c, m$ bounded, 
 $\essinf \Re c > -\frac{(d-2)^2}4$ and $\|m\|_{\infty} \le \varepsilon$ small enough has a unique fundamental solution $\Gamma(t,s)$ defined  on $\L^2_{x,v}$ for $-\infty<s\le t <\infty$.  
 \end{cor}

\begin{proof} This is a perturbation  of the kinetic Schr\"odinger operator in Corollary~\ref{cor:schrodinger}. We can apply  Theorem~\ref{thm:stability} with $\cP$ the multiplication by 
 $m |x|^{-2/3}$ and  
 \[ \Pdot=  \L^2_{t,v}\L^{p,2}_{\vphantom{t^a}x} \cap \L^2_{t,x}\Hdot^1_{\vphantom{t}v}, \quad p=\tfrac{6d}{3d-2},
 \] provided $\|m\|_{\infty}$ is small enough. 
\end{proof}

 \subsection{Non-local diffusions}
 
 Let us now consider some non-local operator $\cA$ when $\beta<1$ satisfying for example
  \eqref{e:ellip-upper-inhom} and \eqref{eq:weakelliptic-bis}.   
  The question is: what kind of lower order perturbation can we add to $\cA$? It is clear from the above examples that we  can look at different ones separately and then together by following the sum/intersection principle described above.

\subsubsection{Addition of a potential $\oc$.} In this case, we may assume  
\[
 \oc \in (\L^{\tilde r}_{t}\L^{\tilde p,\infty}_{x}\L^{\tilde q, \infty}_{v})_{\varepsilon} + \L^\infty_{t,x,v} \ 
\]
  for some $  \tilde r, \tilde p,\tilde q\in (1,\infty]$ satisfying the scaling-critical relation
\[
\frac{2\beta}{\tilde r} +\frac{(2\beta+1)d}{\tilde p} +\frac{d}{\tilde q}= 2\beta
\]
and small enough $\varepsilon>0$. For example,  we may take $\oc\in \L^{(1+1/\beta)d+1}_{t,x,v}$ with arbitrary norm, or $\oc \in \L^\infty_{t,x}\L^{d/2\beta, \infty}_{v}$ with small norm, such as $c|v|^{-2\beta}$. Then we write $\cA+\oc= \widetilde \cA +\cP$ with $\cP$ a perturbation of  $\widetilde \cA$ as in $\mathbf{(H1)}$ and $\cP$ is defined on the 
 space     
\[
\Pdot= \L^{\vphantom{p}r}_{t}\L^{p,2}_{ \vphantom{t}x}\L^{q, 2}_{\vphantom{t} v} \cap \L^2_{t,x}\Hdot_{ v}^\beta,
\]
with $r=2(\tilde r)'$, $p=2(\tilde p)'$, $q=2(\tilde q)'$. Note that
\[
\frac{2\beta}{r} +\frac{(2\beta+1)d}{p} +\frac{d}{q}= (\beta+1)d, \quad r,p,q\in [2,\infty).
\]

The lemma corresponding to Lemma~\ref{lem:ex1} becomes

\begin{lem}\label{lem:ex3}{}\ 
\begin{enumerate}
\item For $p_{0}, q_{0}\in [2,\infty)$ with $\frac{(2\beta+1)d}{p_{0}} +\frac{\vphantom{(2\beta+1)}d}{ q_{0}}=(\beta+1)d-\beta$, 
we have the anisotropic Sobolev embedding
\[
\L^2_{x}\Hdot_{ v}^\beta
\cap  \L^2_{v}\Hdot_{\vphantom{x} x}^{\gb} \subset \L^{p_{0},2}_{\vphantom{t}x}\L^{q_{0},2}_{\vphantom{t}v}. \]
\item For $p,q,r\in [2,\infty)$ with $\frac{2\beta}{r} +\frac{(2\beta+1)d}{p} +\frac{\vphantom{2\beta}d}{q}= (\beta+1)d$,  we have
 the continuous inclusion
 \[(\Wbd)' \subset  \L^{\vphantom{p}r}_{t}\L^{p,2}_{ \vphantom{t}x}\L^{q, 2}_{\vphantom{t} v}.
 \]

\end{enumerate}

\end{lem}

The verification of the perturbative setting under $\mathbf{(H1)}$ is done as in Proposition~\ref{lem:H2}. 

\

If $\tilde r=\tilde p=\infty$, $\tilde q= d/2\beta>1$, and if we assume $\Re \oc \ge 0$, we do not need the smallness of the $\L^\infty_{t,x}\L^{d/2\beta, \infty}_{v}$ norm invoking  $\mathbf{(H2)}$ with space 
\[
\Bd= \L^{\vphantom{p}2}_{t,x}
\L^{q, 2}_{\vphantom{t} v} \cap \L^2_{t,x}\Hbd_{v},  \quad q=\frac{2\beta d}{d-2\beta}.
\] 

\
Another example with same space $\Bd$ comes from Pitt's inequality. 

\begin{cor}[fractional kinetic Schr\"odinger operator]\label{cor:fractionalschrodinger} Assume $\beta<d/2$. The fractional kinetic Schr\"odinger operator 
\[(\partial_{t}+v\cdot\nabla_{x}) +(-\Delta_{v})^{\beta}+ c|v|^{-2\beta}\]
 with $c$ a constant or even a bounded function of $t,x,v$ satisfying
$\essinf \Re c> -1/C_{d,\beta}^2$ with
\[
C_{d,\beta}
=
2^{-\beta}\, 
\frac{
\Gamma\left(\frac{d}{4}-\frac{\beta}{2}\right)
}{
\Gamma\left(\frac{d}{4}+\frac{\beta}{2}\right)
},
\]
 has a unique fundamental solution $\Gamma(t,s)$ defined and contractive on $\L^2_{x,v}$ for $-\infty<s\le t <\infty$.  
 \end{cor}

\begin{proof}
 It follows from the same argument as in  Corollary~\ref{cor:schrodinger} using the Hardy-Rellich  inequality with sharp constant  (the corresponding inequality for Fourier series appears in Pitt \cite{Pitt1937} and the sharp constant independently in the three works \cite{Herbst1977, Beckner1995, Yafaev1999}), namely for 
 $f\in\mathcal{S}(\mathbb{R}^d)$
\[
\left\||v|^{-\beta}f\right\|_{L^2(\mathbb{R}^d)}
\leq
C_{d,\beta}
\left\|(-\Delta_{v})^{\beta/2}f\right\|_{L^2(\mathbb{R}^d)}.
\]
\end{proof}

As in Example 2, we may also exchange the order of the Lorentz norms in the variables $x$ and $v$, allowing for example potentials $\oc$ bounded by  $c|x|^{-2\gb}$ with small constant $c$. 

\

We can also take sums in $(\L^{\tilde r_{1}}_{t}\L^{\tilde p_{1},\infty}_{x}\L^{\tilde q_{1}, \infty}_{v}+ \L^{\tilde r_{2}}_{t}\L^{\tilde q_{2}, \infty}_{v}\L^{\tilde p_{2},\infty}_{x})_{\varepsilon} + \L^\infty_{t,x,v}$, etc. 

\
 
\subsubsection{Fractional perturbations to $\cA$.} We give an example, not claiming exhaustiveness.  Consider 
\[
p_{t,x} (f,g)=\int_{\R^{2d}} \tilde k (t,x,v,w) \frac{f(v)-f(w)}{|v-w|^{\beta+d/2}}\,\overline g(v) \dd (v,w)
\]
and assume $\tilde k\in \L^{\vphantom{\tilde p}2\tilde r}_{t}\L^{2\tilde p, \infty}_{\vphantom{t}x}\L^{\vphantom{\tilde p}\infty}_{\vphantom{t}v}\L^{\vphantom{\tilde p}2}_{\vphantom{t}w}$  with $\tilde r, \tilde p\in (1,\infty]$ and $\frac{2\beta}{\tilde r} +\frac{(2\beta+1)d}{\tilde p} = 2\beta$.
Using Cauchy-Schwarz and the characterization of the $\Hdot^\beta_{v}$ norm (see \cite{DNPV}) we see that 
\[
|p_{t,x} (f,g)|  \lesssim \bigg(\int_{v}|g(v)|^2\|\tilde k(t,x,v,\cdot)\|_{\L^2_{\vphantom{t}w}}^2\dd v\bigg)^{1/2} \|f\|_{\Hdot^\beta_{v}}.
\]
Integrating in the $t,x$ variable and using the Cauchy-Schwarz inequality, 
if $g\in \L^r_{t}\L^{p,2}_{x}\L^2_{v}$, with $r=2(\tilde r)'$, $p=2(\tilde p)'$, we have
\[|\langle \cP f, g\rangle|  \lesssim \|g\|_{\L^{\vphantom{p}r}_{t}\L^{p,2}_{\vphantom{t}x}\L^{\vphantom{p}2}_{\vphantom{t}v}} \|\tilde k\|_{\L^{\vphantom{p} 2\tilde r}_{t}\L^{2\tilde p, \infty}_{\vphantom{t}x}\L^{\vphantom{p}\infty}_{\vphantom{t}v}\L^{\vphantom{p}2}_{\vphantom{t} w}}\|f\|_{\L^{\vphantom{p}2}_{t,x}\Hdot^\beta_{\vphantom{t}v}}.
\]
Note that
\[
\frac{2\beta}{r} +\frac{(2\beta+1)d}{p} +\frac{d}{2}= (\beta+1)d, \quad r,p\in [2,\infty),
\]
so that $(\Wbd)' \subset  \L^{\vphantom{p}r}_{t}\L^{p,2}_{ \vphantom{t}x}\L^{2}_{\vphantom{t} v}
 $ by Lemma~\ref{lem:ex3}, (ii). The space 
  
\[
\Pdot= \L^{\vphantom{p}r}_{t}\L^{p,2}_{ \vphantom{t}x}\L^{2}_{\vphantom{t} v} \cap \L^2_{t,x}\Hdot_{ v}^\beta,
\]
satisfies $\mathbf{(B)}$  and one  can verify $\mathbf{(H1)}$ for $\cA$ and apply the stability theorem to the perturbation  $\cP$ if $\tilde k$ has small norm.

\bibliographystyle{plain}

\begin{thebibliography}{10}

\bibitem{anceschi2026giorginashmosertheorykineticequations}
F.~Anceschi, G.~Palatucci, and M.~Piccinini.
\newblock {D}e {G}iorgi-{N}ash-{M}oser theory for kinetic equations with
  nonlocal diffusions, 2026.
\newblock \url{arXiv:2401.14182}

\bibitem{francesca2021fundamental}
F.~Anceschi and A.~Rebucci.
\newblock {On the fundamental solution for degenerate Kolmogorov equations with
  rough coefficients}.
\newblock {\em Journal of elliptic and parabolic equations}, 9:63--92, 2023.



\bibitem{ACE} W.~Arendt, I.~Chalendar, and  R.~Eymard.
\newblock Lions representation theorem and applications. 
\newblock {\em J. Math. Anal. Appl.}, 522, No. 2, 23 p. (2023). 

\bibitem{Ar67}
D.~G. Aronson. \emph{Bounds for the fundamental solution of a parabolic equation}. Bull. Amer. Math. Soc. \textbf{73} (1967), 890--896.

\bibitem{Ar68}
D.~G. Aronson. \emph{Non-negative solutions of linear parabolic equations}. Ann.
Scuola Norm. Sup. Pisa (3) \textbf{22} (1968), 607--694.

  \bibitem{AB}
P.~Auscher and K.~Baadi.
\newblock Fundamental solutions for parabolic equations and systems: universal
  existence, uniqueness, representation.
\newblock {\em J. Math. Anal. Appl.}, 552(1):Paper No. 129806, 48, 2025.
  

\bibitem{AIN}
P.~Auscher, C.~Imbert, and L.~Niebel.
\newblock Weak solutions to {K}olmogorov-{F}okker-{P}lanck equations:
  regularity, existence and uniqueness, 2024.
\newblock \url{arXiv:2403.17464}

\bibitem{AIN2}
P.~Auscher, C.~Imbert, and L.~Niebel.
\newblock Fundamental solutions to {K}olmogorov-{F}okker-{P}lanck equations with rough coefficients: Existence, uniqueness, upper estimates.
\newblock {\em SIAM Journal on Mathematical Analysis}, 57(2):2114--2137, 2025. 

\bibitem{AN}
P.~Auscher,  and L.~Niebel.
\newblock Kinetic Sobolev spaces, 2026.
\newblock \url{arXiv:2603.17491}

\bibitem{AN2}
P.~Auscher,  and L.~Niebel.
\newblock Existence and upper bounds of fundamental solutions to kinetic equations with rough non-local diffusion, 2026.
\newblock \url{arXiv:2610.06214}

\bibitem{Beckner1995}
W.~Beckner,
\newblock Pitt's inequality and the uncertainty principle,
\newblock \emph{Proceedings of the American Mathematical Society},
123 (1995), no.~6, 1897--1905.







\bibitem{DNPV}
E.~Di~Nezza, G.~Palatucci, and E.~Valdinoci.
\newblock Hitchhiker's guide to the fractional {S}obolev spaces.
\newblock {\em Bull. Sci. Math.}, 136(5):521--573, 2012.


\bibitem{Gaddy}
P.~Gaddy. 
\newblock Nash-Aronson Estimate for the Linear Kinetic Fokker-Planck equation.
\newblock \url{arXiv:2604.22944}




\bibitem{MR3923847}
F.~Golse, C.~Imbert, C.~Mouhot, and A.~F.~Vasseur.
\newblock Harnack inequality for kinetic {F}okker-{P}lanck equations with rough
  coefficients and application to the {L}andau equation.
\newblock {\em Ann. Sc. Norm. Super. Pisa Cl. Sci. (5)}, 19(1):253--295, 2019.



\bibitem{Herbst1977}
I.~W. Herbst,
\newblock Spectral theory of the operator
$(p^2+m^2)^{1/2}-Ze^2/r$,
\newblock \emph{Communications in Mathematical Physics},
53 (1977), 285--294.



\bibitem{Hou2026}
H.~Hou and X.~Zhang.
\newblock Heat kernel estimates for nonlocal kinetic operators.
\newblock {\em Probability Theory and Related Fields}, July 2026.


  
  \bibitem{MR4049224}
C.~Imbert and L.~Silvestre.
\newblock The weak {H}arnack inequality for the {B}oltzmann equation without
  cut-off.
\newblock {\em J. Eur. Math. Soc. (JEMS)}, 22(2):507--592, 2020.





\bibitem{MR3717350}
A.~Lanconelli and A.~Pascucci.
\newblock Nash estimates and upper bounds for non-homogeneous {K}olmogorov
  equations.
\newblock {\em Potential Anal.}, 47(4):461--483, 2017.


\bibitem{Lions}
J.-L.~Lions.
\newblock Sur les probl\`emes mixtes pour certains syst\`emes paraboliques dans
 des ouverts non cylindriques.
\newblock {\em Ann. Inst. Fourier (Grenoble)}, 7:143--182, 1957.
  
  \bibitem{O'Neil1963convolution}
R.~O{'}Neil.
\newblock Convolution operators and {L}(p,q) spaces.
\newblock {\em Duke Math. J.}, 30(1):129--142, 1963.

\bibitem{KN}
K.~Nystr\"om
\newblock The kinetic Kato square root problem.
\newblock \url{arXiv:2611.03104}
	

 \bibitem{PP} 
 {A.~Pascucci,  and A.~Pesce}.
 \newblock Sobolev embeddings for kinetic {Fokker}-{Planck} equations.
 \newblock   {\em J. Funct. Anal.},
  {286},
 ( {7}), 2024, 44 pp.
 
\bibitem{MR2000352}
A.~Pascucci and S.~Polidoro.
\newblock A {G}aussian upper bound for the fundamental solutions of a class of
  ultraparabolic equations.
\newblock {\em J. Math. Anal. Appl.}, 282(1):396--409, 2003.

\bibitem{MR2068847}
A.~Pascucci and S.~Polidoro.
\newblock The {M}oser's iterative method for a class of ultraparabolic
 equations.
\newblock {\em Commun. Contemp. Math.}, 6(3):395--417, 2004.

\bibitem{Pitt1937}
H.~R. Pitt,
\newblock Theorems on Fourier series and power series,
\newblock \emph{Duke Mathematical Journal},
3 (1937), no.~4, 747--755.




\bibitem{SW}
  {E. M.~Stein, and G.~Weiss.} 
\newblock {\em Introduction to Fourier analysis on Euclidean spaces.}
\newblock Princeton Mathematical Series, No. 32. Princeton University Press, Princeton, N.J., 1971.

\bibitem{Tartar98}
L.~Tartar.
\newblock Imbedding theorems of {S}obolev spaces into {L}orentz spaces.
\newblock {\em Boll. Unione Mat. Ital. Sez. B Artic. Ric. Mat. (8)},
  1(3):479--500, 1998.

\bibitem{Yafaev1999}
D.~R. Yafaev,
\newblock Sharp constants in the Hardy--Rellich inequalities,
\newblock \emph{Journal of Functional Analysis},
168 (1999), no.~1, 121--144.

\end{thebibliography}

\end{document}